\documentclass [12pt, twoside, reqno] {amsart}

\usepackage{texdraw}
\usepackage{amsfonts}
\usepackage{amssymb}
\usepackage{latexsym}
\usepackage{epsf,graphicx}
\usepackage{epsf,graphicx,latexsym,%
}

\newcommand{\pl}{\partial}
\newcommand{\pr}{\noindent{\bf Proof.}\quad }
\newcommand{\epr}{\ $\blacksquare$ \vspace{3mm} }

\newcommand{\be} {\begin{eqnarray}}
\newcommand{\ee} {\end{eqnarray}}
\newcommand{\bep} {\begin{eqnarray*}}
\newcommand{\eep} {\end{eqnarray*}}

\newcommand {\Hol}{\mathop{\rm Hol}\nolimits}
\renewcommand {\Im}{\mathop{\rm Im}\nolimits}
\renewcommand {\Re}{\mathop{\rm Re}\nolimits}

\newcommand {\GP}{\mathcal{G}_{\alpha,\beta}(\Pi)}
\newcommand {\GD}{\mathcal{G}_{\alpha,\beta}(\Delta)}

\newcommand{\C}{{\mathbb C}}

\newtheorem{remar}{Remark}[section]
\newtheorem{examp}{Example}[section]
\newtheorem{defin}{Definition}[section]
\newtheorem{corol}{Corollary}[section]
\newtheorem{propo}{Proposition}[section]
\newtheorem{theorem}{Theorem}[section]
\newtheorem{lemma}{Lemma}[section]

\newtheorem{conj}{Conjecture}

\newcommand{\rema}{\begin{remar}\rm}
\newcommand{\erema}{$\blacktriangleright$\end{remar}}

\newcommand{\exa}{\begin{examp}\rm}
\newcommand{\eexa}{$\blacktriangleright$\end{examp}}

\def\lwvec(#1 #2){\linewd 0.1
           \lvec(#1 #2)
           \linewd 0.05}

\title[Parabolic type semigroups]
    {Parabolic type semigroups: asymptotics and order of contact}

\author[M. Elin]{Mark Elin}

\address{Department of Mathematics,
         Ort Braude College,
         Karmiel 21982,
         Israel}

\email{mark$\_$elin@braude.ac.il}

\author[F. Jacobzon]{Fiana Jacobzon}

\address{Department of Mathematics,
         Ort Braude College,
         Karmiel 21982,
         Israel}

\email{fiana@braude.ac.il}

\begin{document}

\maketitle

\begin{abstract}

We study the asymptotic behavior of parabolic type semigroups
acting on the unit disk as well as those acting on the right
half-plane. We use the asymptotic behavior to investigate the
local geometry of the semigroup trajectories near the boundary
Denjoy--Wolff point. The geometric content includes, in
particular, the asymptotes to trajectories, the so-called limit
curvature, the order of contact, and so on. We then establish
asymptotic rigidity properties for a broad class of semigroups of
parabolic type.

\vspace{4mm}

{\footnotesize Key words and phrases: holomorphic mapping,
asymptotic behavior, parabolic type semigroup, contact order,
rigidity.

}
\end{abstract}

\section{Preliminaries}

The theory of semigroups of holomorphic self-mappings of a given
domain in the complex plane $\mathbb{C}$ has been developed
intensively over the last few decades. The study began with the
basic work of E. Berkson and H. Porta \cite{B-P} (see, e.g.,
\cite{SD} and \cite{E-S-book} for a recent state of this theory).
This paper is devoted to the study of a wide class of parabolic
type semigroups acting on the open unit disk and on the right
half-plane.

Throughout the paper, $\Hol(D,\C)$ denotes the set of holomorphic
functions on a domain $D\subset\mathbb{C}$ and $\Hol(D)$ denotes
the set of holomorphic self-mappings of $D$. Recall that \textit{a
one-parameter continuous semigroup} (semigroup, for short) acting
on $D$ is  a family $S=\left\{F_t\right\}_{t\geq 0}\subset\Hol(D)$
such that
\begin{flushleft}\begin{itemize}
\item[(i)] $F_{t}(F_{s}(z))= F_{t+s}(z)$ for all $t,s\geq 0$ and
$z\in D$,

\item[(ii)] $\lim\limits_{t\rightarrow 0^+}F_{t}(z)=z$ for all
$z\in D$.
\end{itemize}
\end{flushleft}
Berkson and Porta \cite{B-P} proved that  each semigroup acting on
$D$ when $D$ is either the open unit disk $\Delta=\{
z\in\mathbb{C}:|z|<1\}$ or the right half-plane $\Pi=\{
z\in\mathbb{C}: \Re z>0\}$ is differentiable with respect to
$t\in\mathbb{R}^{+}=[0,\infty )$. Thus, for each one-parameter
continuous semigroup the limit
\[
\lim_{t\rightarrow 0^{+}}\frac{F_{t}(z)-z}{t}=f(z),\quad z\in D,
\]
exists and defines a holomorphic function
$f\in\Hol(D,\mathbb{C})$. This function $f$ is called the
\textit{(infinitesimal) generator of} $S$. Moreover, the function
$ u(t,z):=F_{t}(z), \, (t,z)\in\mathbb{R}^{+}\times D$, is the
unique solution of the Cauchy problem
\[
\left\{
\begin{array}{l}
{\displaystyle\frac{\partial u(t,z)}{\partial t}} = f(u(t,z))
,\vspace{3mm}\\ u(0,z)=z,\quad z\in D.
\end{array}
\right.
\]

In the same paper, Berkson and Porta proved that
$f\in\Hol(\Delta,\C)$ is a semigroup generator if and only if
there exist a point $\tau\in\overline{\Delta}$ and a function
$p\in\Hol(\Delta,\C)$ with $\Re p(z)\ge0$, such that
\[
f(z)=(\tau-z)(1-z\bar{\tau })p(z).
\]
This representation is unique. Moreover, if $S$ contains neither
the identity mapping nor an elliptic automorphism of $\Delta$,
then $\tau$ is a unique attractive fixed point of $S$, i.e.,
$\lim\limits_{t\to\infty}F_t(z)=\tau$ for all $z\in\Delta,$ and
$\lim\limits_{r\to1^-}F_t(r\tau)=\tau$. The point $\tau$ is called
the \textit{Denjoy--Wolff point of} $S$.

Recently the asymptotic behavior of semigroups including the local
geometry of semigroup trajectories near their boundary
Denjoy--Wolff point $\tau\in\partial\Delta$ has attracted
considerable attention. It was shown in \cite{E-S1} that if
$\tau\in\pl\Delta$, then the angular derivative of $f$ at
$\tau\in\pl\Delta$ defined by $\displaystyle f'(\tau)=
\angle\lim\limits_{z\to\tau}\frac{f(z)}{z-\tau}$ exists and is a
non-positive real number.

There is an essential difference between semigroups whose
generator $f$ satisfies $f'(\tau)<0$ (semigroups of hyperbolic
type) and those whose generator $f$ satisfies $f'(\tau)=0$
(semigroups of parabolic type). For example, in the hyperbolic
case, the rate of convergence of the semigroup to its
Denjoy--Wolff point is exponential, while in the parabolic case,
the convergence is slower. The main problem we address can be
stated as follows.

{\it Determine the rate of convergence of parabolic type
semigroups; more precisely, find the asymptotic expansion up to a
term small enough}.

\vspace{2mm}

Obviously, every semigroup trajectory $\gamma_z=\left\{F_t(z),\
t\ge0 \right\},\ z\in\Delta,$ is an analytic curve. Thus, the
tangent line and the circle of curvature at each its point
$F_s(z)$ exist and move as $s$ increases. The following natural
question arises.

{\it Do tangent lines and disks of curvature have, in some sense,
a limit location as $s\to\infty$?}

\vspace{2mm}

It turns out that in the hyperbolic case, limit tangent lines
always exist and depend on the initial point of the trajectory. On
the other hand, for all studied classes of parabolic type
semigroups, all trajectories have the same limit tangent line, but
even its existence has not been proven in general. More precisely,
M. D. Contreras and S. D\'\i az-Madrigal in \cite{C-DM} considered
the set $\mathrm{Slope}^+(\gamma_z)$ of accumulation points (as
$t\to\infty$) of the function $t\mapsto \arg\left(1-\bar\tau
F_t(z)\right)$ and proved that these sets do not depend on
$z\in\Delta$. There are cases in which it is known that
$\mathrm{Slope}^+(\gamma_z)$  is a singleton. The question as to
whether, in general, $\mathrm{Slope}^+(\gamma_z)$ is a singleton
is still open (see \cite{C-DM, E-R-S-Y, E-S-Y, E-K-R-S, E-S2011}
for details).

To be more concrete, we henceforth assume without losing any
generality, that $\tau=1$. We mention (see \cite{E-R-S-Y}) that if
the generator $f$ of a parabolic type semigroup
$S=\{F_t\}_{t\ge0}$ admits the representation
\[
f(z)=a(1-z)^2+o((1-z)^2),
\]
then for each $z\in\Delta$, the limit tangent line to the
trajectory $\gamma_z =\{F_t(z),\ t\ge0\}$ exists, and
\[
\lim_{t\to\infty}\arg(1-F_t(z))=-\arg a.
\]
Hence, this limit depends on neither $z\in\Delta$ nor the
remainder $o((z-1)^2)$. This fact was generalized in \cite{E-S-Y}
(see also \cite{E-K-R-S}) for the case
\[
f(z)=a(1-z)^{1+\alpha}+o((1-z)^{1+\alpha})\quad\mbox{for some
}\alpha>0.
\]
Moreover, it was shown in \cite{E-S-Y} that $$\alpha\le 2,\quad
\left\vert \arg a\right\vert \leq \frac{\pi }{2}\min
\{\alpha,2-\alpha\},\quad
\lim_{t\rightarrow \infty }\arg (1-F_{t}(z))=-\frac{1}{\alpha}
\arg a.
$$
In particular, this implies that all the trajectories are tangent
to the unit circle if and only if $\alpha\le1$ and $\arg a =\pm
\frac{\pi\alpha}{2}$ (see \cite{E-K-R-S} for more details).
Proposition~\ref{th_asy_one_term} below completes these results.
An advanced question in this study is the following.

{\it How close is a semigroup trajectory to its tangent line?}

\vspace{2mm}

Following \cite{E-S2011}, for each $z\in\Delta$, we denote the
curvature of the trajectory $\gamma_z$ at the point $F_s(z)$  by
$\kappa(z,s)$ and define the \textit{limit curvature of the
trajectory} by ${\kappa(z):=\lim\limits_{s\to\infty}\kappa(z,s)}$,
if the limit exists. Therefore, the above question can be reduced
to the following one.

{\it When is the limit curvature finite?}

\vspace{2mm}

This question was studied in \cite{E-S2011}, where it was shown
that every trajectory of a hyperbolic type semigroup has a finite
limit curvature, while the finiteness in the parabolic case is, in
a sense, exceptional. Namely, it was proved in \cite{E-S2011} that
if a semigroup generator is $(3+\varepsilon)$-smooth at the
Denjoy--Wolff point, in the sense that it admits the
representation
\[
f(z)=a(1-z)^2+ b(1-z)^3+R(z),
\]
where $R\in\Hol(\Delta,\mathbb{C}),\ \displaystyle
\lim\limits_{z\to1}\frac{R(z)}{(1-z)^{3+\varepsilon}}=0$, and
$a\not=0$, {\it

(a) if $\displaystyle\ \Im\frac{b}{a^2}\not=0$, then all of the
trajectories 
have infinite limit curvature, i.e., $\kappa(z)=\infty$ for all
$z\in\Delta$;

(b) if $\displaystyle\ \Im\frac{b}{a^2}=0$, the limit curvature of
every trajectory $\gamma_z$ is finite. The value $\kappa(z)$ was
calculated explicitly in \cite{E-S2011}}.

Thus, under the above assumptions, if $\kappa(z)$ is finite for
some $z\in\Delta$, then it must be finite for all $z\in\Delta$.

Once again, we see that there is a cardinal difference between
semigroups of hyperbolic and parabolic types. In the hyperbolic
case under some smoothness conditions, the limit curvature is
always finite; in the parabolic case, the limit curvature may be
infinite. For the above reasons, for parabolic type semigroups, a
more relevant question is

{\it find the contact order of a trajectory and the limit tangent
line} (which is less than $2$ when the limit curvature is
infinite).

This problem leads to the so-called rigidity problem, which is
that of finding the weakest conditions on two holomorphic mappings
at a boundary point under which the mappings coincide. Beginning
with the outstanding work of D.~Burns and S.~G.~Krantz \cite{B-K},
this problem has attracted considerable interest (see \cite{SD-06,
E-L-R-S} and reference therein). As a rule, the rigidity problem
for one-parameter semigroups is approached by looking for
conditions on generators. Another approach is related to semigroup
asymptotics. Here, we investigate the rigidity problem via contact
order of the trajectories. In our setting, the next question is
natural.

{\it What is the minimal   contact order of trajectories of
parabolic type semigroups required to ensure that the semigroups
coincide?}

\vspace{2mm}

We solve the above problems for parabolic type semigroups
$S=\{F_t\}_{t\ge0}$ whose generators $f\in\Hol(\Delta,\C)$ admit
the representation
\begin{equation}\label{circ_1a-e}
f(z)=a(1-z)^{1+\alpha}+R(z),
\end{equation}
or the representation
\begin{equation}\label{circ_4a-e}
f(z)=a(1-z)^{1+\alpha}+b(1-z)^{1+\alpha+\beta}+R_1(z),
\end{equation}
where $\alpha\in(0,2],\ \beta>0,\ a\not=0$, and functions
$R,R_1\in\Hol(\Delta,\mathbb{C})$ satisfy
\begin{equation}\label{remainders}
\lim\limits_{z\to1}\frac{R(z)}{(1-z)^{1+\alpha}}=0,
\quad\mbox{and}\quad
\lim\limits_{z\to1}\frac{R_1(z)}{(1-z)^{1+\alpha+\beta}}=0.
\end{equation}
As previously mentioned, if $|\arg a|< \frac{\pi\alpha}{2}$, the
semigroup converges non-tangentially. Since, we use formulas
\eqref{circ_1a-e} and \eqref{circ_4a-e} to expand $f(F_t(z))$, in
the case $|\arg a|\not = \frac{\pi\alpha}{2}\,,$ the limits in
\eqref{remainders} can be replaced by angular limits.

In what follows, $\GD$ denotes the set of semigroup generators
having the form \eqref{circ_4a-e} with $a\not=0$ and function
$R_1$ satisfying \eqref{remainders}.

Also, we apply a linearization model given by Abel's functional
equation
\begin{equation}\label{abel}
h\left(F_t(z)\right)=h(z)+t.
\end{equation}
It is rather easy to see that the function $h:\Delta\mapsto\C$
defined by
\begin{equation}\label{h}
h'(z)f(z)= 1,\quad h(0)=0,
\end{equation}
solves functional equation (\ref{abel}). This function is
univalent and, due to (\ref{abel}), is convex in the positive
direction of the real axis . Sometimes $h$ is called the {\it
K{\oe}nigs function} for the semigroup (see \cite{C-DM, E-R-S-Y,
E-S-Y, Sis} and \cite{E-S-book}).

The class of semigroups acting on $\Pi$ and the class acting on
$\Delta$ are conjugated by $\Phi_t(w)=C\circ F_t\circ C^{-1}(w),$
where $C$ is the Cayley transform $C(z)=\frac{1+z}{1-z}$. For
technical reasons, we first study the behavior of semigroups
acting on $\Pi$. Whence $S=\{F_t\}_{t\geq 0}$ has Denjoy--Wolff
point $\tau=1$, semigroup $\Sigma=\left\{\Phi_t\right\}_{t\geq
0}\subset\Hol(\Pi)$ has Denjoy--Wolff point $\infty$, its
generator $\phi$ belongs to $\Hol(\Pi,\overline{\Pi})$, and the
semigroup $\Sigma=\left\{\Phi_t\right\}_{t\geq 0}$ satisfies the
Cauchy problem
\[
\left\{
\begin{array}{l}
{\displaystyle\frac{\partial \Phi_t(w)}{\partial
t}}=\phi\left(\Phi_t(w)\right),\vspace{3mm}\\
\left.\Phi_t(w)\right|_{t=0}=w,\quad w\in \Pi.
\end{array}
\right.
\]

We modify the K{\oe}nigs function $h$ defined by (\ref{h}) to
$\sigma:=h\circ C^{-1}$. Direct calculations show that for all
$w\in\Pi$, this modified function satisfies Abel's functional
equation
\begin{equation}\label{abel1}
\sigma\left(\Phi_t(w)\right)=\sigma(w)+t
\end{equation}
as well as the initial value problem
\begin{equation}\label{h-1}
\sigma'(w)\phi(w)=1, \quad \sigma(1)=0.
\end{equation}

It follows by from Julia's Lemma (see, for example, \cite{SJH-93,
SD, E-S-book}) that since the Denjoy-Wolff point of $\Sigma$ is
$\infty$, hence $\Re \Phi_t(w)$ is an increasing function in $t$
for $t\ge0$. This prompts an additional question.

{\it What conditions ensure the existence of asymptotes to
semigroup trajectories?}

\vspace{2mm}

Note in passing that {\it a semigroup trajectory  $\gamma_z
\subset \Delta$ has a finite limit curvature if and only if
$C(\gamma_z) \subset \Pi$ has an asymptote as $t\to\infty$.}

\vspace{3mm}

In Section~\ref{sect_hp}, we study the asymptotic behavior of
semigroups acting on $\Pi$. These semigroups not only give us a
machinery for our main results, but are of intrinsic interest.
Despite the fact that these semigroups tend to $\infty$, the
asymptotic behavior which we describe enables us to distinguish
those semigroups whose trajectories are either {\it asymptotically
parallel, or mutually convergent, or mutually divergent} (see
Definition~\ref{def_traj} below). As an application, we deduce the
rather surprising result that in the case $\alpha< \min \{1,\beta
\}$, the motion on each trajectory is accelerating. Consequently,
the distance between two particles starting at different points of
the same trajectory grows to $\infty$ (see Corollary
\ref{col_conv_div}). In addition, we present a complete
description of conditions for the existence of asymptotes to
semigroup trajectories and their possible coincidence.

In Section~\ref{sect_disk},  we turn to semigroups acting on
$\Delta$ generated by functions of the class $\GD$.
Theorem~\ref{th_main_asy} contains a full description of the
asymptotic behavior of such semigroups. One of the phenomena
discovered is that for semigroup generators of the
form~\eqref{circ_4a-e}, if $\beta\le\alpha$, the first two terms
of the asymptotic expansion of the generated semigroup do not
depend on the initial point. On the other hand, if $\beta>\alpha$,
the initial point affects as from the second term. Moreover, if
$\beta\le\alpha$, it may happen that all the trajectories have the
same contact order (see Definition~\ref{def_order} below), while
if $\beta>\alpha$, there exists a trajectory $\gamma$ of maximal
contact order. One of the geometric implications of this
phenomenon is that there exists no semigroup trajectory lying
between $\gamma$ and the limit tangent line. Each trajectory
starting from a point between $\gamma$ and the tangent line must
intersect the tangent and approach it from the side opposite from
$\gamma$ (see Remark~\ref{rem1} below). We also provide conditions
under which the limit curvature is either zero, finite, or
infinite.

In Section~\ref{sect_rigid}, we study the contact order of two
trajectories and use results from earlier sections to establish
rigidity criteria for parabolic type semigroups. As a bonus, we
discover another interesting geometric phenomenon. In the case
$0<\beta <\alpha$, each trajectory is closer to all other
trajectories than to their common limit tangent line. Thus, all
the trajectories approach this tangent line from the same side
(see Remark \ref{rem2}).

\bigskip

\section{Semigroups on the right half-plane}\label{sect_hp}

\setcounter{equation}{0}

In this section, we study parabolic type semigroups acting on the
right half-plane $\Pi$.  We begin by assuming that only the first
term in the asymptotic expansion of the generator is known.

\begin{lemma}\label{th-cur-par01}
Let $\{\Phi_{t}\}_{t\geq 0}\in\Hol(\Pi)$ be a semigroup of
parabolic type with the Denjoy--Wolff point at $\infty$ generated
by mapping $\phi$. Suppose that
\begin{equation}\label{asy_gen}
\phi(w)=A(w+1)^{1-\alpha}+\varrho(w),
\end{equation}
where $\varrho\in\Hol(\Pi,\mathbb{C})$, and
$\displaystyle\angle\lim_{w\to\infty}\frac{\varrho(w)}{(w+1)^{1-\alpha}}=0.
$ Then
\begin{equation}\label{arg}
\Phi_t(w)=(\lambda
t)^{\frac{1}{\alpha}}+\Gamma(w,t)\quad\mbox{with }
\lim_{t\to\infty}t^{-\frac{1}{\alpha}}\Gamma(w,t)=0,
\end{equation}
and
\begin{equation}\label{difference1}
\lim_{t\to\infty}
\left(\Phi_t^{\alpha}(w)-\Phi_t^{\alpha}(1)\right)=\lambda\sigma(w),
\end{equation}
where $\sigma$ is defined by (\ref{h-1}) and $\lambda=\alpha A.$\\
\end{lemma}
As already mentioned, if a semigroup generator satisfies
\eqref{asy_gen} then $0<\alpha\leq 2$. The case $\alpha=1$ was
considered in \cite[Theorem 4.1(i)]{E-S2011}.

\pr Fix $w\in\Pi$ and consider $\Phi_t(w)$ as a (complex valued)
function of the real variable~$t$. Since
$\lim_{t\to\infty}\Phi_t(w)=\infty$, L'H\^{o}pital's rule gives
\begin{eqnarray} \label{lim1}
\lim_{t\to\infty}\frac{\left(\Phi_t(w)+1\right)^{\alpha}}{t+1}=
\lim_{t\to\infty}\frac{\alpha\left(\Phi_t(w)+1\right)^{\alpha-1}\phi\left(\Phi_t(w)\right)}1
\nonumber\\
\lim_{t\to\infty}\alpha\left(
A+\frac{\varrho\left(\Phi_t(w)\right)}{\Phi_t^{1-\alpha}(w)}\right)
=\lambda.
\end{eqnarray}
Thus,
\[
\lim_{t\to\infty}\frac{\Phi_t(w)}{t^{\frac{1}{\alpha}}}=\lim_{t\to\infty}\frac{\Phi_t(w)}{t^{\frac{1}{\alpha}}}\cdot\frac{\Phi_t(w)+1}{\Phi_t(w)}\cdot\left(\frac{t}{t+1}\right)^{\frac{1}{\alpha}}=\lim_{t\to\infty}\frac{\left(\Phi_t(w)+1\right)}{(t+1)^{\frac{1}{\alpha}}}=\lambda^{\frac{1}{\alpha}}.
\]
This proves (\ref{arg}). Furthermore,
\begin{eqnarray*}
&&\lim_{t\to\infty}
\left(\Phi_t^{\alpha}(w)-\Phi_t^{\alpha}(1)\right) =
\lim_{t\to\infty}\int_1^w\left(\Phi_t^{\alpha}(z)\right)'dz
\\
&=&\lim_{t\to\infty}\int_1^w \alpha\Phi_t^{\alpha-1}(z)
\frac{\phi\left(\Phi_t(z) \right)}{\phi(z)}dz \\
&=&\lim_{t\to\infty}\int_1^w \left(\frac{\Phi_t(z)
}{\Phi_t(z)+1}\right)^{\alpha-1}\frac{\alpha}{\phi(z)}\left(A+\frac{\varrho\left(\Phi_t(z)\right)}{\Phi_t^{1-\alpha}(z)}
\right)dz \\
&=& \lambda\int_1^w \frac{dz}{\phi(z)}=\lambda\sigma(w)
\end{eqnarray*}
by (\ref{lim1}). \epr

In the case in which the function $\rho$ in (\ref{asy_gen}) can be
written as $\rho(w)=B(w+1)^{1-\alpha-\beta}+\varrho_1(w)$ with
$\beta>0$ and $\lim\limits_{w\to\infty}
\frac{\varrho_{1}(w)}{(w+1)^{1-\alpha-\beta}}=0,$ we can obtain a
more precise estimate for the asymptotic behavior of the generated
semigroup. Denote the set of generators $\phi \in \Hol (\Pi,
\overline{\Pi})$ of the form
\begin{equation}\label{repr_cp10less}
\phi(w)=A(w+1)^{1-\alpha}+B(w+1)^{1-\alpha-\beta}+\varrho_1(w)
\end{equation}
by $\GP$, where $\varrho_{1}\in\Hol(\Pi,\C)$ satisfies
${\lim\limits_{w\to\infty} \frac{\varrho_{1}(w)}
{(w+1)^{1-\alpha-\beta}}=0}$ and ${A \not=0}$. For the remainder
of this section, we deal with semigroups whose Denjoy--Wolff point
is $\tau=\infty$ and whose infinitesimal generators lie in $\GP$.
We also set
\begin{equation}\label{lambda_mu}
\lambda=\alpha A\quad\mbox{and}\quad \mu=\frac{B}{A}\,.
\end{equation}

It turns out that semigroups have different asymptotic behavior
depending on whether $\beta<\alpha$, $\beta=\alpha$, or
$\beta>\alpha$. We start with the case $\beta<\alpha$.

\begin{theorem}\label{th-cur-par03}
Let $\Sigma=\{\Phi_{t}\}_{t\geq 0}\in\Hol(\Pi)$ be a semigroup
generated by a mapping $\phi \in \GP$ with $0<\beta<\alpha\leq 2$.
Then
\begin{equation}\label{repr_cp_new3less}
\Phi_t(w)+1=(\lambda t)^{\frac{1}{\alpha}}\left(1+\frac{\mu}
{\alpha-\beta}(\lambda t)^{-\frac{\beta}{\alpha}}
+\Gamma(w,t)\right),
\end{equation}
where
$
\lim\limits_{t\to\infty}\displaystyle
t^{\frac{\beta}{\alpha}}\Gamma(w,t) =0.
$
\end{theorem}

\pr First we show that
\begin{equation}\label{limit15less}
\lim_{t\to\infty}\frac{1}{(t+1)^{1-\frac{\beta}{\alpha}}} \left(
\left(\Phi_t(w)+1\right)^{\alpha}-\lambda t-\frac{\alpha
\mu}{\alpha-\beta}(\lambda(t+1))^{1-\frac{\beta}{\alpha}}\right)=0.
\end{equation}
Using (\ref{repr_cp10less}), we calculate
\begin{eqnarray*}
&&\frac{d}{ds}\left( \left(\Phi_s(w)+1\right)^{\alpha}-\lambda s-\frac{\alpha \mu}{\alpha-\beta}(\lambda(s+1))^{1-\frac{\beta}{\alpha}}\right)\\
\\
&=&\alpha\left(B \left( (\Phi_s(w)+1)^{-\beta}
-(\lambda(s+1))^{-\frac{\beta}{\alpha}}\right)
+\left(\Phi_s(w)+1\right)^{\alpha-1}\rho_{1}(\Phi_s(w))\right)
\end{eqnarray*}
\begin{equation}\label{der_calculation}
=\frac{\alpha}{(s+1)^{\frac{\beta}{\alpha}}}
\left[B\left(\mu(s,w)-\lambda^{-\frac{\beta}{\alpha}}\right)+
\frac{\mu(s,w)
\rho_{1}(\Phi_s(w))}{\left(\Phi_s(w)+1\right)^{1-\alpha-\beta}}\right],
\end{equation}
where $\mu(s,w)=
\left(\frac{(s+1)^{\frac{1}{\alpha}}}{\Phi_s(w)+1}\right)^{\beta}$.
Also, by (\ref{lim1}),
$\lim\limits_{s\to\infty}\mu(s,w)=\lambda^{-\frac\beta\alpha}$. By
our assumption on $\rho_1$ in (\ref{repr_cp10less}), it follows
that
\[
\lim_{s\to\infty}
\left[B\left(\mu(s,w)-\lambda^{-\frac{\beta}{\alpha}}\right) +
\frac{\mu(s,w)\rho_{1}(\Phi_s(w))}{\left(\Phi_s(w)+1\right)^{1-\alpha-\beta}}\right]=0.
\]
Therefore, for each $\varepsilon>0$, there exists $t_0$ such that
for all $s>t_0$,
\[
\left|B\left(\mu(s,w)-\lambda^{-\frac{\beta}{\alpha}}\right)+
\frac{\mu(s,w)\rho_{1}(\Phi_s(w))}{\left(\Phi_s(w)+1\right)^{1-\alpha-\beta}}\right|
<\frac{\varepsilon (\alpha-\beta)}{\alpha^2 }\,,
\]
while, for each $0\leq s \leq t_{0}$, there exists $K>0$ such that
\[
\left|B\left(\mu(s,w)-\lambda^{-\frac{\beta}{\alpha}}\right)+
\frac{\mu(s,w)\rho_{1}(\Phi_s(w))}{\left(\Phi_s(w)+1\right)^{1-\alpha-\beta}}\right|
<\frac{K (\alpha-\beta)}{\alpha^2 }\,.
\]
Since
\begin{eqnarray}\label{integral}
&&\left(\Phi_t(w)+1\right)^{\alpha}-\lambda t-\frac{\alpha
\mu}{\alpha-\beta}(\lambda(t+1))^{1-\frac{\beta}{\alpha}}
-(w+1)^{\alpha}+\frac{\alpha
\mu}{\alpha-\beta}\lambda^{1-\frac{\beta}{\alpha}} \nonumber \\
&&=\int_0^{t}\frac{d}{ds}\left(
\left(\Phi_s(w)+1\right)^{\alpha}-\lambda s-\frac{\alpha
\mu}{\alpha-\beta}(\lambda(s+1))^{1-\frac{\beta}{\alpha}}\right)ds,
\end{eqnarray}
formula (\ref{der_calculation}) implies
\begin{eqnarray*}
&& \frac{1}{(t+1)^{1-\frac{\beta}{\alpha}}}
\left|\left(\Phi_t(w)+1\right)^{\alpha}-\lambda t-\frac{\alpha
\mu}{\alpha-\beta}(\lambda(t+1))^{1-\frac{\beta}{\alpha}} \right|
\\  && \leq \frac{1}{(t+1)^{1-\frac{\beta}{\alpha}}}
\left[\int_{0}^{t_0}\frac{K(\alpha-\beta)}{\alpha(s+1)^{\frac{\beta}{\alpha}}}
ds +\int_{t_0}^{t}\frac{\varepsilon
(\alpha-\beta)}{\alpha(s+1)^{\frac{\beta}{\alpha}}}ds\right.\\
&& \left.+\left|(w+1)^{\alpha}-\frac{\alpha
\mu}{\alpha-\beta}\lambda^{1-\frac{\beta}{\alpha}}\right|\right] \\
&&
=\varepsilon+(K-\varepsilon)\frac{(t_0+1)^{1-\frac{\beta}{\alpha}}}{(t+1)^{1-\frac{\beta}{\alpha}}}+\frac
{K+\left|(w+1)^{\alpha}-\frac{\alpha
\mu}{\alpha-\beta}\lambda^{1-\frac{\beta}{\alpha}}\right|}{(t+1)^{1-\frac{\beta}{\alpha}}}.
\end{eqnarray*}
Since $\varepsilon>0$ is arbitrary, (\ref{limit15less}) follows
and
\begin{equation}\label{r3less}
\left(\Phi_t(w)+1\right)^{\alpha}=\lambda t+\frac{\alpha
\mu}{\alpha-\beta}
(\lambda(t+1))^{1-\frac{\beta}{\alpha}}+\Gamma_1(w,t)
\end{equation}
with $\lim\limits_{t\to\infty}\displaystyle
(t+1)^{\frac{\beta}{\alpha}-1}\Gamma_1(w,t)=0.$

To proceed we calculate
\begin{eqnarray*}
\Phi_t(w)+1&=&\Bigl(\lambda t+\frac{\alpha
\mu}{\alpha-\beta}(\lambda(t+1))^{\frac{\alpha-\beta}{\alpha}}+\Gamma_1(w,t)\Bigr)^{\frac{1}{\alpha}}
\\
&=&(\lambda t)^{\frac{1}{\alpha}}\left(1+\frac{\alpha
\mu}{\alpha-\beta}\frac{1}{\lambda t}
(\lambda(t+1))^{\frac{\alpha-\beta}{\alpha}}+\frac{\Gamma_1(w,t)}{\lambda
t}\right)^{\frac{1}{\alpha}}
\\
&=&(\lambda
t)^{\frac{1}{\alpha}}\left(1+\frac{\mu}{\alpha-\beta}\frac{1}{\lambda
t}
(\lambda(t+1))^{\frac{\alpha-\beta}{\alpha}}+{\widetilde{\Gamma}_1(w,t)}\right),
\end{eqnarray*}
where $\lim\limits_{t\to\infty}\displaystyle
 t^{\frac{\beta}{\alpha}}\widetilde{\Gamma}_1(w,t)=0.$ This proves the assertion. \epr

\vspace{3mm}

The case $\beta=\alpha$ can be treated similarly. We state the
analogous result.

\begin{theorem}[cf., Theorem 4.1(ii) in \cite{E-S2011}]\label{th-cur-par02}
Let $\Sigma=\{\Phi_{t}\}_{t\geq 0}\subset \Hol(\Pi)$ be a
semigroup generated by a mapping $\phi \in \GP$ with
$0<\beta=\alpha \le 2$, i.e.,
\[
\phi(w)=A(w+1)^{1-\alpha}+B(w+1)^{1-2\alpha}+\varrho_1(w), \quad
\lim\limits_{w\to\infty}
\displaystyle\frac{\varrho_{1}(w)}{(w+1)^{1-2\alpha}}=0.
\]
Then
\begin{equation}\label{repr_cp_new1}
\Phi_t(w)+1=\left(\lambda
t\right)^{\frac{1}{\alpha}}\left(1+\frac{\mu}{\alpha} \cdot
\frac{\log(t+1)}{\lambda t}+\Gamma(w,t)\right),
\end{equation}
where $\lim\limits_{t\to\infty}\displaystyle\frac{t \Gamma(w,t)}{
\log(t+1)}=0$.
\end{theorem}

It turns out that the asymptotic behavior of a semigroup can be
estimated more precisely when the remainder $\varrho_1$
in~(\ref{repr_cp10less}) satisfies a stronger condition. The next
result generalizes \cite[Theorem 4.1(iii)]{E-S2011}.

\begin{propo}\label{addit_asymp}
Let $\phi \in \GP$ be given by \eqref{repr_cp10less}, where
$\varrho_{1}\in\Hol(\Pi,\C)$ satisfies $\lim\limits_{w\to\infty}
\displaystyle (w+1)^{2\alpha-1+\varepsilon} \varrho_{1}(w)=0$ for
some positive $\varepsilon$, and let $\Sigma=\{\Phi_{t}\}_{t\geq
0}\subset\Hol(\Pi)$ be a semigroup generated by $\phi$.
\begin{itemize}
\item[(i)] If $\beta=\alpha$, then there exists a constant $C$
such that for all $w \in \Pi$,
\[
\left(\Phi_t(w)+1\right)^{\alpha}=\lambda t+\mu \log(t+1)+\lambda
\sigma(w)+ C+{\Gamma}(w,t),
\]
where $\lim\limits_{t \to 0} {\Gamma}(w,t)=0$.

\item[(ii)] If $\frac\alpha 2<\beta<\alpha$ then there exists a
constant $C$ such that for all $w \in \Pi$,
\[
\left(\Phi_t(w)+1\right)^{\alpha}=\lambda t+\frac{\mu
\alpha}{\alpha-\beta} (\lambda t)^{1-\frac{\beta}{\alpha}}+\lambda
\sigma(w)+ C+{\Gamma}(w,t),
\]
where $\lim\limits_{t\to\infty}{\Gamma}(w,t)=0$.
\end{itemize}
\end{propo}

\pr Since the proofs of assertions (i) and (ii) are similar, we
prove only assertion (ii). We first show that the limit
\[
H(w):=\lim_{t\to\infty}\left(\left(\Phi_t(w)+1\right)^{\alpha}-\lambda
t-\frac{\alpha
\mu}{\alpha-\beta}(\lambda(t+1))^{1-\frac{\beta}{\alpha}}\right)
\]
exists for each $w \in \Pi$. Indeed, by the calculations in
(\ref{der_calculation}) and (\ref{integral}), we have
\begin{eqnarray*}
&&\left(\Phi_t(w)+1\right)^{\alpha}-\lambda t-\frac{\alpha
\mu}{\alpha-\beta}(\lambda(t+1))^{1-\frac{\beta}{\alpha}}
=\\
& =&\int_{0}^{t}\alpha B \left( (\Phi_s(w)+1)^{-\beta}
-(\lambda(s+1))^{-\frac{\beta}{\alpha}}\right)ds \\
& +&
\int_{0}^{t}\alpha\left(\Phi_s(w)+1\right)^{\alpha-1}\rho_{1}(\Phi_s(w))ds+(w+1)^{\alpha}-\frac{\alpha
\mu}{\alpha-\beta}(\lambda(t+1))^{1-\frac{\beta}{\alpha}}.
\end{eqnarray*}
For the first integral, (\ref{r3less}) implies
\begin{eqnarray*}
\frac{1}{(\Phi_s(w)+1)^{\beta}}
-\frac{1}{(\lambda(s+1))^{\frac{\beta}{\alpha}}}
=\frac{1}{(\lambda(s+1))^{\frac{\beta}{\alpha}}}
\left(\left(\frac{\lambda(s+1)}{(\Phi_s(w)+1)^{\alpha}}\right)^{\frac{\beta}{\alpha}}-1
\right) \\
=\frac{1}{(\lambda(s+1))^{\frac{\beta}{\alpha}}}\left(\left(1+
\frac{\alpha\mu}{\alpha-\beta}(\lambda(s+1))^{-\frac\beta\alpha}
-\frac{1}{s+1}+\frac{\Gamma_1(w,s)}{\lambda(s+1)}\right)^{-\frac{\beta}{\alpha}}-1
\right).
\end{eqnarray*}

Since this expression is
$O\left(\displaystyle\frac{1}{(s+1)^{\frac{2
\beta}{\alpha}}}\right)$ with $ \displaystyle\frac{2
\beta}{\alpha}>1$, the first integral converges. For the second
integral, we have
\begin{eqnarray*}
&&\frac{\alpha\rho_{1}(\Phi_s(w))}{\left(\Phi_s(w)+1\right)^{1-\alpha}}
= \alpha\left(\Phi_s(w)+1\right)^{-\alpha-\varepsilon}
\frac{\rho_{1}(\Phi_s(w))}{\left(\Phi_s(w)+1\right)^{1-2\alpha-\varepsilon}} \\
&&=\frac{\alpha}{(s+1)^{\frac{\alpha+\varepsilon}{\alpha}}}
\cdot\left(\frac{s+1}{\left(\Phi_s(w)+1\right)^{\alpha}}\right)^{\frac{\alpha+\varepsilon}{\alpha}}
\cdot\frac{\rho_{1}(\Phi_s(w))}{\left(\Phi_s(w)+1\right)^{1-2\alpha-\varepsilon}}
\, .
\end{eqnarray*}
If
$\displaystyle\lim_{t\to\infty}\frac{\rho_{1}(\Phi_s(w))}{\left(\Phi_s(w)+1\right)^{1-2\alpha-\varepsilon}}=0$,
then by (\ref{lim1}) we have
\[
\alpha
\left|\frac{\rho_{1}(\Phi_s(w))}{\left(\Phi_s(w)+1\right)^{1-\alpha}}\right|=O
\left(\frac{1}{(s+1)^{1+\frac{\varepsilon}{\alpha}}}\right),
\]
and so the second integral also converges. The proof of the
equality
\[H(w)-H(1)=\lim_{t\to\infty}
\left(\left(\Phi_t(w)+1\right)^{\alpha}-\left(\Phi_t(1)+1\right)^{\alpha}\right)=\lambda\sigma(w)
\]
is similar to that of (\ref{difference1}). This proves that
\[
\left(\Phi_t(w)+1\right)^{\alpha}=\lambda t+\frac{\mu
\alpha}{\alpha-\beta}
(\lambda(t+1))^{1-\frac{\beta}{\alpha}}+\lambda \sigma(w)+
C+\Gamma(w,t),
\]
where $C=H(1)$ and
\[
\Gamma(w,t)=\left(\Phi_t(w)+1\right)^{\alpha}-\lambda t-\frac{\mu
\alpha}{\alpha-\beta}
(\lambda(t+1))^{1-\frac{\beta}{\alpha}}-\lambda \sigma(w)-C.
\]
The result now follows, since
$\lim\limits_{t\to\infty}\Gamma(w,t)=0$. \epr

We now turn to the case $\beta>\alpha$.

\begin{theorem}\label{th-cur-par04}
Let $\Sigma=\{\Phi_{t}\}_{t\geq 0}\in\Hol(\Pi)$ be a semigroup
generated by $\phi \in \GP$, where $\beta\not=\alpha$ satisfies
$0<k\alpha\le\beta<(k+1)\alpha$ for some $k\in\mathbb{N}$. Then
\begin{eqnarray}\label{repr_cp_new3}
&&\Phi_t(w)+1 =(\lambda t)^{\frac{1}{\alpha}} \times \\
&&\ \times \nonumber
\left[1+ \sum_{j=1}^k\left(\begin{array}c 1/{\alpha} \\
j\end{array}\right) \left(\frac{\sigma_1(w)}t\right)^j
 + \frac{\mu}
{\alpha-\beta}(\lambda t)^{-\frac{\beta}\alpha} +
\Gamma(w,t)\right],
\end{eqnarray}
where $ \lim\limits_{t\to\infty}\displaystyle
t^{\frac{\beta}{\alpha}} \Gamma(w,t) =0$ and
$\sigma_1(w)=\sigma(w)+\frac{2^\alpha}{\lambda} -
\int\limits_1^\infty\left(\sigma'(v)-
\frac{1}{A}(v+1)^{\alpha-1}\right)dv$.
\end{theorem}

\pr We apply the K\oe nigs function $\sigma$ which satisfies
(\ref{abel1}). By (\ref{h-1}),
\[
\sigma'(w)=\frac{1}{\phi(w)}=\frac{1}{A(w+1)^{1-\alpha}+B(w+1)^{1-\alpha-\beta}+\varrho_1(w)}.
\]
A direct calculation gives
\begin{equation}\label{sigma_prime}
\sigma'(w)=\frac{1}{A}(w+1)^{\alpha-1}-\frac{B}{A^2}(w+1)^{\alpha-1-\beta}+r(w),
\end{equation}
where ${\lim\limits_{w\to\infty} r(w){(w+1)^{1-\alpha+\beta}}=0}$.
In particular, this implies that the improper integral
$\int\limits_1^\infty r(v)dv$ converges, say to $C_1$. An
application of L'H\^{o}pital's rule gives
\[
\lim\limits_{w\to\infty} \frac{\int_1^w
r(v)dv-C_1}{(w+1)^{\alpha-\beta}}=0.
\]
It now follows from (\ref{sigma_prime}) that
\[
\sigma(w)=\int_1^w\sigma'(v)dv=\frac{1}{\lambda}(w+1)^\alpha
+\frac{\alpha \mu }{\lambda(\beta-\alpha)}(w+1)^{\alpha-\beta}+
C_2 + r_1(w),
\]
where ${\lim\limits_{w\to\infty} r_1(w){(w+1)^{\beta-\alpha}}=0}$
and $C_2=\int\limits_1^\infty\left(\sigma'(v)-
\frac{1}{A}(v+1)^{\alpha-1}\right)dv -\frac{2^\alpha}{\lambda} $.
Substituting this asymptotic expansion into Abel's functional
equation (\ref{abel1}) yields
\begin{equation}\label{mid_asymp}
\frac{1}{\lambda}(\Phi_t(w)+1)^\alpha +\frac{\alpha \mu
}{\lambda(\beta-\alpha)}(\Phi_t(w)+1)^{\alpha-\beta}+ C_2 +
r_1(\Phi_t(w)) = \sigma(w)+t.
\end{equation}
According to Lemma~\ref{th-cur-par01},
$\lim\limits_{t\to\infty}\frac{(\Phi_t(w)+1)^{\alpha-\beta}}{(t+1)^{1-\frac\beta\alpha}}
=\lambda^{1-\frac\beta\alpha}$. Therefore,
\[
(\Phi_t(w)+1)^{\alpha-\beta} =(\lambda(t+1))^{1-\frac\beta\alpha}
+ r_2(t,w),
\]
where $\lim\limits_{t\to\infty}(t+1)^{\frac\beta\alpha-1}
r_2(t,w)=0$. Comparing the last relation with (\ref{mid_asymp}),
we conclude that
\begin{eqnarray*}
(\Phi_t(w)+1)^\alpha &=& \lambda t + \lambda(\sigma(w)-C_2)+
\frac{\alpha \mu
} {\alpha-\beta}(\lambda(t+1))^{1-\frac\beta\alpha}  \\
&&-\,\frac{\alpha\mu}{\beta-\alpha}r_2(t,w)
-\lambda r_1(\Phi_t(w)) \\
&=& \lambda t + \lambda\sigma_1(w)+\frac{\alpha \mu
}{\alpha-\beta}(\lambda(t+1))^{1-\frac\beta\alpha} +R (t,w),
\end{eqnarray*}
where $\lim\limits_{t\to\infty}R(t,w){(t+1)^{\frac\beta\alpha -1}}
=0$ and $\sigma_1(w)=\sigma(w)-C_2 $.

The formula for the sum of a binomial series gives
\begin{eqnarray*}
&&\Phi_t(w)+1=\Bigl(\lambda t + \lambda\sigma_1(w)+\frac{\alpha
\mu }{\alpha-\beta}(\lambda(t+1))^{1-\frac\beta\alpha} +R
(t,w)\Bigr)^{\frac{1}{\alpha}}
\\
&=&(\lambda t)^{\frac{1}{\alpha}} \left(1+\frac{\sigma_1(w)}t
\right)^{\frac{1}{\alpha}} \left(1+ \frac{\alpha\mu }
{\lambda(\alpha-\beta)}\frac{(\lambda(t+1))^{1-\frac\beta\alpha}}{t+\sigma_1(w)}
+ \frac{R (t,w)}{t+\sigma_1(w)}\right)^{\frac{1}{\alpha}}
\\
&=&(\lambda t)^{\frac{1}{\alpha}}\left(1+
\sum_{j=1}^k\left(\begin{array}c 1/{\alpha} \\ j\end{array}\right)
\left(\frac{\sigma_1(w)}t \right)^j
+O\left(\frac1{t^{k+1}}\right)\right)\times \\
&&\hspace{1cm} \times\left(1+ \frac{\mu } {\lambda(\alpha-\beta)}
\frac{(\lambda(t+1))^{1-\frac\beta\alpha}}{t+\sigma_1(w)} +
o\left(t^{-\frac\beta\alpha}\right) \right)\\
& =& (\lambda t)^{\frac{1}{\alpha}}\left(1+
\sum_{j=1}^k\left(\begin{array}c 1/{\alpha} \\ j\end{array}\right)
\left(\frac{\sigma_1(w)}t \right)^j + \frac{\mu (\lambda
t)^{-\frac\beta\alpha} } {\alpha-\beta} +
o\left(t^{-\frac\beta\alpha}\right)\right).
\end{eqnarray*}
The proof is complete. \epr

Theorems \ref{th-cur-par03}--\ref{th-cur-par04} \ give more than
asymptotic expansions of semigroups. Using standard methods of
analysis we can deduce, on the basis of these theorems,
interesting facts about the geometry of semigroup trajectories.
For example, we give criteria on $\alpha$ and $\beta$ which ensure
the existence/non-existence of asymptotes to semigroup
trajectories. We also determine whether the asymptote exists for
all initial points $w\in\Pi$ or only for $w$ from some subset of
$\Pi$, and whether the asymptote (if it exists) depends on the
initial point. As we will see below, the cases in which the
asymptote passes through $-1$ are of special interest.

First, we decompose the set $\Omega=\{(\alpha,\beta): \ 0<\alpha
\le 2, \ \beta>0  \}$ of all possible pairs of the parameters into
the following subsets:
\begin{eqnarray*}
&&\Omega_1:= \left \{(\alpha,\beta) \in \Omega : \, \alpha,\beta>1
\right \},\\
&&\Omega_2:= \left \{(\alpha,\beta) \in \Omega : \, 1=\alpha <
\beta \right \},\\
&&\Omega_3:= \left \{(\alpha,\beta) \in \Omega : \, \alpha < \min
\{1,\beta\}\right \},\\
&&\Omega_4:= \left \{(\alpha,\beta) \in \Omega : \, 1= \beta
<\alpha \leq 2\right\},\\
&&\Omega_5:= \left \{(\alpha,\beta) \in \Omega : \, \beta \leq
\min \{1,\alpha\}  \right \} \setminus \Omega_4.
\end{eqnarray*}
Obviously, these sets are pairwise disjoint and their union covers
$\Omega$ (see Fig. 1).

\begin{figure}\centering

    \includegraphics[angle=0,width=7cm,totalheight=4cm]{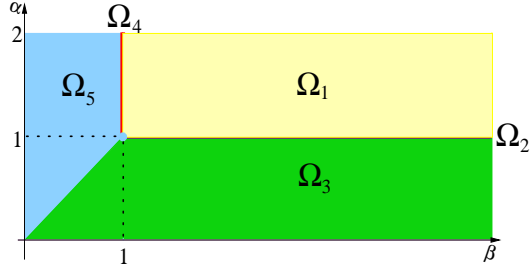}
    \caption{Partition of the set of parameters}
\end{figure}

\begin{propo}\label{cor_asymptotes}
Let $\Sigma=\{\Phi_{t}\}_{t\geq 0}\subset\Hol(\Pi)$ be a semigroup
generated by $\phi \in \GP$.
\begin{itemize}
\item[(i)] If  $(\alpha,\beta) \in \Omega_1$, then all the
trajectories of $\Sigma$ have the same asymptote. This asymptote
passes through the point $-1$.

\item[(ii )] If $(\alpha,\beta) \in \Omega_2$, then each
trajectory has its own asymptote. The asymptote depends on the
initial point.

\item[(iii)] If $(\alpha,\beta) \in \Omega_3$, then the only
trajectory $\gamma$ defined by the condition $\left.\Im \sigma_1
\right|_\gamma =0$ has an asymptote. This asymptote passes through
the point $-1$.

\item[(iv)] If  $(\alpha,\beta) \in \Omega_4$, then all the
trajectories of $\Sigma$ have the same asymptote. This asymptote
passes through the point $-1$ if and only if $\Im
\left(BA^{-\frac{\alpha+\beta}{\alpha}} \right)=0$.

\item[(v)] If  $(\alpha,\beta) \in \Omega_5$ and $\Im
\left(BA^{-\frac{\alpha+\beta}{\alpha}} \right) =0$, then all the
trajectories have the same asymptote. This asymptote passes
through the point $-1$.

\item[(vi)] If  $(\alpha,\beta) \in \Omega_5$ and $\Im
\left(BA^{-\frac{\alpha+\beta}{\alpha}} \right) \neq 0$, no
trajectory of $\Sigma$ has an asymptote.
\end{itemize}
\end{propo}

\pr The problem reduces to an examination of the limit
\[
\lim\limits_{t\to\infty}\Im\left( \overline{\lambda^\frac1\alpha}
(\Phi_t(w)+1) \right).
\]
Indeed, the trajectory $\{\Phi_t(w)\}_{t\geq 0}$ has an asymptote
if and only if this limit exists finitely. Moreover, if this limit
vanishes, the asymptote passes through the point $-1$. We
determine the existence of this limit and its value (if it it
exists) using asymptotic expansions \eqref{repr_cp_new3less},
\eqref{repr_cp_new1} and \eqref{repr_cp_new3}.

For the case $\beta<\alpha$, using formula
\eqref{repr_cp_new3less} from Theorem \ref{th-cur-par03} we obtain
\[
\lim\limits_{t\to\infty}\Im\left( \overline{\lambda^\frac1\alpha}
(\Phi_t(w)+1)
\right)=\lim\limits_{t\to\infty}\frac{|\lambda|^\frac 2 \alpha
t^\frac {1-\beta} {\alpha}}{\alpha-\beta}\Im \left(\mu
\lambda^{-\frac \beta \alpha}\right).
\]
The limit on the right vanishes for all pairs $(\alpha,\beta)$
with $1<\beta<\alpha$; hence the asymptote exists and passes
through $-1$. If $1=\beta<\alpha$, the same limit exists and the
asymptote passes through the point $\frac{|\lambda|^\frac 2 \alpha
}{\alpha-1}\Im\left( \mu \lambda^{-\frac 1 \alpha}\right)-1$.
Finally, if $\beta<\min \{1,\alpha \}$, then an asymptote exists
if and only if $\Im \mu \lambda^{-\frac \beta \alpha}=0$, and if
it does, it passes through $-1$. This proves assertion (iv) and
parts of assertions (i) and (v). The remaining parts of assertions
(i) and (v) as well as assertions (ii), (iii) and (iv) follow from
a similar argument using Theorems
\ref{th-cur-par02}--\ref{th-cur-par04}. \epr

The particular case of assertion (v) for $\alpha=\beta=1$ was
treated in \cite[Theorem 4.2(a),(b)]{E-S2011}.

Another interesting issue is to estimate how far are two
trajectories of the same semigroup having different initial
points. The theorems above immediately imply the following.

\begin{corol}\label{cor_1}
Let $\Sigma=\{\Phi_{t}\}_{t\geq 0}\in\Hol(\Pi)$ be a semigroup
generated by $\phi \in \GP$. For all $w\in\Pi$,

\begin{itemize}

\item[(i)] if $\beta=\alpha$, then $\lim\limits_{t\to\infty}
\frac{t^{1-\frac{1}{\alpha}}}{\log(t+1)}\left(\Phi_t(w) -\Phi_t(1)
\right)=0$;

\item[(ii)] if $\beta<\alpha$, then $\lim\limits_{t\to\infty}
t^{\frac{\beta-1}{\alpha}}\left(\Phi_t(w)-\Phi_t(1)\right)=0$;

\item[(iii)] if $\beta>\alpha$, then $\lim\limits_{t\to\infty}
\left(t^{\frac{\beta-1}{\alpha}}\left(\Phi_t(w)-\Phi_t(1)\right)-
\frac{1}{\alpha}\lambda^\frac 1 \alpha t^{\frac {\beta}{\alpha}-1}
\sigma(w)\right)=0$.

\end{itemize}
\end{corol}

In turn, Corollary~\ref{cor_1} yields a simple description of the
relative position of the semigroup trajectories going to the
Denjoy--Wolff point at infinity. To formulate it we introduce the
following notions.

\begin{defin}\label{def_traj}
Let $\Sigma=\{\Phi_t\}_{t\ge0}\subset\Hol(\Pi)$ be a semigroup
with the Denjoy--Wolff point at infinity. For $w_1,w_2 \in \Pi$,
let
\[
s(w_1,w_2):=\lim\limits_{t\to\infty}
\left(\Phi_t(w_1)-\Phi_t(w_2)\right)
\]
whenever the limit exists. We say that the semigroup trajectories
are
\begin{itemize}
\item[(i)]  mutually convergent if $s\equiv 0$ on $\Pi\times\Pi$,

\item[(ii)] asymptotically parallel if $s$ is well defined on $\Pi
\times \Pi$ and does not vanish on
$\Pi\times\Pi\setminus\{w_1=w_2\}$,

\item[(iii)] mutually divergent if $s(w_1,w_2)=\infty$ for all
$w_1\neq w_2$.

\end{itemize}
\end{defin}

Note that if the trajectories are mutually convergent, then for
every compact $K \subset \Pi$ and $\varepsilon>0$ there exists
$t_0$ such that for each $t>t_0$, the set $\left\{\Phi_t(w), w \in
K \right \}$ is contained in a disk of radius $\varepsilon$.

\begin{corol}[see Fig. 2]\label{col_conv_div}
Under conditions of Corollary~\ref{cor_1}, the following
assertions hold.

\begin{itemize}
\item[(i)] If $\alpha>1$ and $\beta\ge1$, then the trajectories of
$\Sigma$ are mutually convergent.

\item[(ii)] If $\beta>\alpha=1$, then all the trajectories of
$\Sigma$ are asymptotically parallel. Moreover, the function
$\frac{s(w_1,w_2)}{\sigma(w_1)-\sigma(w_2)}$ is constant.

\item[(iii)] If $\alpha< \min \{1,\beta \}$, then all the
trajectories of $\Sigma$ are mutually divergent. In particular,
\[
\lim\limits_{t\to\infty}
\left(\Phi_{t+1}(w)-\Phi_t(w)\right)=\infty \quad \mbox {for all}
\quad w \in \Pi.
\]
\end{itemize}
\end{corol}

\begin{figure}\centering

    \includegraphics[angle=0,width=9cm,totalheight=4cm]{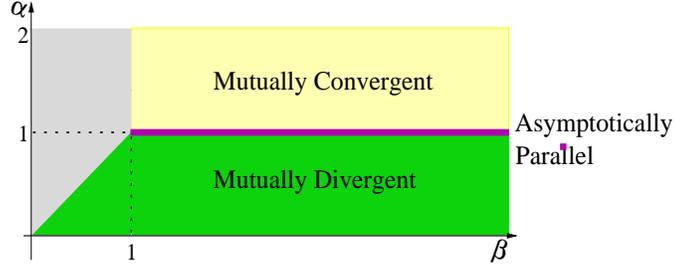}
    \caption{Relative position of trajectories}
\end{figure}

\bigskip

\section{Semigroups on the unit disk}\label{sect_disk}
\setcounter{equation}{0}

In this section, we  study the asymptotic behavior of parabolic
type semigroups. The conclusions derived in \cite{E-S-Y} and
\cite{E-S2011} are a specific case of the results below, which can
be applied to a broader set of semigroups.

Let $S$ be a semigroup of holomorphic self-mappings of the open
unit disk $\Delta$. Using the Cayley transform $\displaystyle
C(z)=\frac{1+z}{1-z}, $\, we transfer the study of semigroups
acting on $\Delta$ to that of those acting on $\Pi$. For a given
semigroup $S=\left\{F_t\right\}_{t\geq 0}\subset\Hol(\Delta)$ with
Denjoy--Wolff point $\tau=1$, we construct the semigroup
$\Sigma=\left\{\Phi_t\right\}_{t\geq 0}\subset\Hol(\Pi)$ with
Denjoy--Wolff point $\infty$ by the composition
\begin{equation}\label{Phi-a}
\Phi_t(w)=C\circ F_t\circ C^{-1}(w).
\end{equation}

Clearly, $\Phi_t \circ C(z) =C \circ F_t(z)$, and hence
\begin{equation}\label{hp-disk}
\Phi_t(C(z))+1= \frac{2}{1 - F_t(z)}.
\end{equation}

If $S$ is continuous (hence, differentiable) in $t$, then so is
$\Sigma$. Suppose that $p\in\Hol(\Delta,\overline{\Pi})$ and
$f(z)=(1-z)^2p(z),\ z \in \Delta,$ generates $S$. Differentiating
$\Phi_t$ given by (\ref{Phi-a}) at $t=0^+$, we conclude that
$\Sigma$ is generated by the mapping
\begin{equation}\label{gen2}
\phi(w)=2p \left(C^{-1}(w)\right).
\end{equation}

Suppose, in addition, that $f$ is of the form~(\ref{circ_1a-e}) or
(\ref{circ_4a-e}). By formula~\eqref{gen2}, the function $\phi$
can be represented, respectively, by
\begin{equation}\label{phi1a}
\phi(w)= 2^{\alpha}a(w+1)^{1-\alpha} +\rho(w)
\end{equation}
with
$\displaystyle\lim\limits_{w\to\infty}\frac{\rho(w)}{(w+1)^{1-\alpha}}=0$,
or by
\begin{equation}\label{phi1b}
\phi(w)= 2^{\alpha}a(w+1)^{1-\alpha} +
2^{\alpha+\beta}b(w+1)^{1-\alpha-\beta}+\rho_1(w)
\end{equation}
with
$\displaystyle\lim\limits_{w\to\infty}\frac{\rho_1(w)}{(w+1)^{1-\alpha-\beta}}=0$.
Thus, $\phi$ has the form (\ref{asy_gen}) or
(\ref{repr_cp10less}), respectively, with
\begin{equation}\label{A-B}
A=2^{\alpha}a\quad\mbox{ and }\quad B=2^{\alpha+\beta}b.
\end{equation}
We also use \eqref{A-B} and modify \eqref{lambda_mu} to
$\lambda=2^\alpha\alpha a$ and $\mu=2^\beta\frac{b}{a}$.

The Cayley transform allows us to apply the results of the
previous section for semigroups acting on $\Pi$ to study
semigroups acting on $\Delta$. The next result is a generalization
of Theorem 1.4(i) in \cite{E-S2011}.

\begin{propo}\label{th_asy_one_term}
Let $S=\left\{F_t\right\}_{t\ge0}$ be a semigroup of holomorphic
self-mappings of $\Delta$ generated by
\begin{equation}\label{circ_1a}
f(z)=a(1-z)^{1+\alpha}+R(z),
\end{equation}
where $R\in\Hol(\Delta,\mathbb{C}),\ \displaystyle
\lim\limits_{z\to1}\frac{R(z)}{(1-z)^{1+\alpha}}=0$. Then
\begin{equation}\label{repr_sh}
\frac1{1-F_t(z)}=\frac{1}{2}(\lambda t)^{\frac{1}{\alpha}}  +
r(z,t)\quad \mbox{with}\quad \lim\limits_{t\to\infty}\displaystyle
t^{-\frac1{\alpha}} r(z,t) =0
\end{equation}
and
\[
\lim_{t\to\infty} \left(\frac{1}{(1-F_t
(z))^\alpha}-\frac{1}{(1-F_t (0))^\alpha}\right)=\frac{\lambda
h(z)}{2^\alpha},
\]
where $h$ is the K\oe nigs function defined by (\ref{h}).
\end{propo}

\pr Substituting (\ref{hp-disk}) into formula (\ref{arg}) yields
\[
\frac1{1-F_t(z)}=\frac{1}{2}(\lambda t)^{\frac{1}{\alpha}}  +
r(z,t),
\]
with  $r(z,t)=\frac 1 2\Gamma(C^{-1} (z),t)$  and
$\lim\limits_{t\to\infty}\displaystyle t^{-\frac1{\alpha}} r(z,t)
=0. $

As in the proof of (\ref{difference1}),
\begin{eqnarray*}
&&\lim_{t\to\infty}
\left((\Phi_t(w)+1)^{\alpha}-(\Phi_t(1)+1)^{\alpha}\right) =
\lim_{t\to\infty}\int_1^w\left((\Phi_t(z)+1)^{\alpha}\right)'dz
\\
&=&\lim_{t\to\infty}\int_1^w
\frac{\alpha}{\phi(z)}\left(A+\frac{\varrho\left(\Phi_t(z)\right)}{\Phi_t^{1-\alpha}(z)}
\right)dz = \lambda\int_1^w \frac{dz}{\phi(z)}=\lambda\sigma(w)
\end{eqnarray*}
with $\sigma(w)=h(C^{-1}(w))$. Again from formula~(\ref{hp-disk}),
we conclude that
\[
\lim_{t\to\infty} \left(\frac{1}{(1-F_t
(z))^\alpha}-\frac{1}{(1-F_t (0))^\alpha}\right)=\frac{\lambda
h(z)}{2^\alpha}.
\]\epr

Using Theorems \ref{th-cur-par03}--\ref{th-cur-par04}, we can
deduce the following asymptotic representation of parabolic type
semigroups for all possible pairs $(\alpha,\beta) \in\Omega$ (for
$\alpha=\beta=1$, cf.,  assertion (ii) of Theorem 1.4 in
\cite{E-S2011}).

\begin{theorem}\label{th_main_asy}
Let $S=\left\{F_t\right\}_{t\ge0}$ be a semigroup of holomorphic
self-mappings of $\Delta$ and let $f \in \GD$ be its generator.
\begin{itemize}
\item[(i)] If $0<\beta<\alpha$, then
\begin{equation}\label{as_less}
\frac1{1-F_t(z)}=\frac{1}{2}(\lambda t)^{\frac{1}{\alpha}}
\left(1+\frac{\mu} {\alpha-\beta}(\lambda
t)^{-\frac{\beta}{\alpha}} + r_1(z,t)\right)
\end{equation}
with $\lim\limits_{t\to\infty}\displaystyle
t^{\frac{\beta}{\alpha}} r_1(z,t) =0$.

\item[(ii)] If $\beta=\alpha$, then
\[
\frac1{1-F_t(z)}=\frac{1}{2}(\lambda t)^{\frac{1}{\alpha}}
\left(1+\frac{\mu}{\alpha} \cdot \frac{\log(t+1)}{\lambda t}+
r_1(z,t)\right)
\]
with $\lim\limits_{t\to\infty}\displaystyle\frac{r_1(z,t)t}
{\log(t+1)}=0$.

\item[(iii)] If $\beta\not=\alpha$ and
$k\alpha\le\beta<(k+1)\alpha$ for some $k\in\mathbb{N}$, then
\begin{eqnarray*}
\frac1{1-F_t(z)}=\frac{1}{2}(\lambda t)^{\frac{1}{\alpha}}
\left[1+ \sum_{j=1}^k\left(\begin{array}c 1/{\alpha} \\
j\end{array}\right) \left(\frac{h_1(z)}t\right)^j \right.\nonumber
\\
\left. + \frac{\mu} {\alpha-\beta}(\lambda
t)^{-\frac{\beta}\alpha} + r_1(z,t)\right]
\end{eqnarray*}
with $\lim\limits_{t\to\infty}\displaystyle
t^{\frac{\beta}{\alpha}} r_1(z,t) =0$ and
$$h_1(z)=h(z)+\frac{2^\alpha}{\lambda} -
\int\limits_0^1\left(h'(s) -
\frac{1}{a(1-s)^{\alpha+1}}\right)ds.$$
\end{itemize}
\end{theorem}

\pr To prove these assertions,  we use Theorems
\ref{th-cur-par03}--\ref{th-cur-par04}. Substituting formulas
(\ref{hp-disk}) and (\ref{A-B}) into~(\ref{repr_cp_new3less}),\
(\ref{repr_cp_new1}) and (\ref{repr_cp_new3}) gives assertions
(i), (ii) and (iii), respectively. To complete the proof, we only
note that (\ref{A-B}) implies that $\lambda=\alpha A$ and
$\mu=\frac{B}{A}=\frac{2^\beta b}{a}$\,, and that the relation for
the K\oe nigs functions $\sigma(w)=h\left(\frac{w-1}{w+1}\right)$
implies $\int\limits_1^\infty\left(\sigma'(v)-
\frac{1}{A}(v+1)^{\alpha-1}\right)dv =\int\limits_0^1\left(h'(s) -
\frac{1}{a(1-s)^{\alpha+1}}\right)ds$. \epr

Theorem \ref{th_main_asy} not only provides a specification of the
asymptotic behavior of semigroups, but also enables us to study
the local geometry of semigroup trajectories in more detail. As
already mentioned, Proposition~\ref{th_asy_one_term} (see also
\cite{E-S-Y}) implies that all trajectories of a semigroup
generated by a function of the form~(\ref{circ_1a}) have the same
limit tangent line. A more detailed analysis requires the
following notion.

\begin{defin}\label{def_order}
Let $\gamma,\gamma^*:[0,\infty)\mapsto\Delta$ be smooth disjoint
curves which satisfy
$\lim\limits_{t\to\infty}\gamma(t)=\lim\limits_{t\to\infty}\gamma^*(t)=
1$. Denote by $d(t)$ the distance between $\gamma(t)$ and
$\gamma^*$. We say that the contact order between $\gamma$ and
$\gamma^*$ (at the point $1$) is $\kappa\ (\kappa\ge0)$, if the
limit
\[
\lim_{t\to\infty}\frac{d(t)}{|1-\gamma(t)|^{1+\kappa}}
\]
exists finitely and is different from zero. If this limit is zero,
we say that the contact order is greater than $\kappa$. In the
case $\gamma^*$ is the limit tangent line of $\gamma$, instead of
``contact order between $\gamma$ and $\gamma^*$" we say ``contact
order of $\gamma$".
\end{defin}

Note that the existence of the limit tangent line guarantees that
the contact order is greater than zero, while the contact order of
a curve $\gamma$ is equal to or greater than $1$ if and only if
$\gamma$ has a finite limit curvature.

\begin{theorem}\label{th_order}
Let $S=\left\{F_t\right\}_{t\ge0}$ be a semigroup of holomorphic
self-mappings of $\Delta$ whose generator $f$ is in $\GD$.

\begin{itemize}
\item[(i)] If $0<\beta<\alpha$, then the contact order of all the
trajectories is at least $\beta$. In the case $\Im
\left(\lambda^{-\frac{\beta}{\alpha}}\mu\right)\not=0$, this order
equals $\beta$.

\item[(ii)] If $\beta=\alpha$, then
\[
\lim_{t\to\infty}\frac{d(t)}{\left|1-F_t(z)
\right|^{1+\alpha}\log|1-F_t(z)|}=
\left|\frac{\lambda^{\frac1\alpha}}{2a}\right|\cdot \left|\Im
\left( \frac{\mu}{\alpha\lambda} \right) \right|;
\]
hence, for any $\varepsilon>0$, the contact order of all the
trajectories is greater than $\alpha-\varepsilon$.

\item[(iii)] If $\beta>\alpha$, then for each $z\in\Delta$ such
that $\Im h_1(z)\not=0$, the trajectory passing through $z$ has
contact order $\alpha$.

For the trajectory $\gamma$ defined by $\Im h_1|_\gamma=0$, the
contact order is at least $\beta$. In particular, if $\Im
\left(\lambda^{-\frac{\beta}{\alpha}}\mu\right)\not=0$, the
contact order is $\beta$.
\end{itemize}
\end{theorem}

\pr By Proposition~\ref{th_asy_one_term},
\[
F_t(z)=1-\frac1{\frac 1 2 \left(\lambda t\right)^{\frac1\alpha} +
r(z,t)}\quad\mbox{with}\quad \lim\limits_{t\to\infty}\displaystyle
t^{-\frac1{\alpha}} r(z,t) =0
\]
and $\lim\limits_{t\to\infty} t^{\frac1\alpha} \left(1-F_t(z)
\right)=2 \lambda^{-\frac1\alpha}$  (see also \cite{E-S-Y}).
Therefore, all the trajectories have the common limit tangent line
$\displaystyle\ell=\left\{z=1-\frac {2x}{\lambda^{\frac1\alpha}},\
x\in \mathbb{R} \right\}$. Following Definition~\ref{def_order},
given a point $z\in \Delta$, we denote the distance between
$F_t(z)$ and $\ell$ by $d(t)$ . Standard analysis yields
\[
d(t)=\frac{\left|\Im\left( \lambda^{\frac1\alpha}\overline{r(z,t)}
\right)\right| }{\left|\lambda^{\frac1\alpha}\right| \left|\frac 1
2 \left(\lambda t\right)^{\frac1\alpha} + r(z,t) \right|^2} ,
\]
so that
\begin{equation}\label{dista}
d(t)= \left|\lambda^{\frac1\alpha}\right|\cdot
\left|\Im\left(\frac{r(z,t) }{\lambda^{\frac1\alpha}}
\right)\right| \cdot \left|1-F_t(z) \right|^2.
\end{equation}

Let $0<\beta<\alpha$. Theorem~\ref{th_main_asy}~(i)  implies that
\[
r(z,t)=\frac{(\lambda t)^{\frac{1-\beta}{\alpha}}}2 \left(
\frac{\mu} {\alpha-\beta} + r_1(z,t)\right),
\]
where $\lim\limits_{t\to\infty} r_1(z,t) =0.$ Hence
\[
d(t)= \frac 1 2 \left|\lambda^{\frac1\alpha}\right|\cdot
\left|\Im\left(\lambda^{-\frac{\beta}{\alpha}}
 \left( \frac{\mu} {\alpha-\beta} +
r_1(z,t)\right) \right)\right| \cdot
t^{\frac{1-\beta}{\alpha}}\left|1-F_t(z) \right|^2
\]
and
\[
\lim_{t\to\infty}\frac{d(t)}{\left|1-F_t(z) \right|^{1+\beta}}=
\frac 1 2\left|\lambda^{\frac1\alpha}\right|\cdot \left|\Im
\frac{\lambda^{-\frac{\beta}{\alpha}}\mu} {\alpha-\beta} \right|.
\]
Assertion (i) follows.

Let $\beta=\alpha$. According to Theorem~\ref{th_main_asy}~(ii),
\[
r(z,t)= \frac 1 2\lambda^{\frac{1}{\alpha}}
t^{\frac{1}{\alpha}-1}\log(t+1) \left( \frac{\mu}{\alpha\lambda} +
r_1(z,t)\right)
\]
with $\lim\limits_{t\to\infty} r_1(z,t) =0$. By formula
(\ref{dista}),
\[
d(t)= \frac12 \left|\lambda^{\frac1\alpha}\right|\cdot \left|\Im
\left( \frac{\mu}{\alpha\lambda} + r_1(z,t)\right) \right| \cdot
t^{\frac{1}{\alpha}-1}\log(t+1) \left|1-F_t(z) \right|^2.
\]
Furthermore,
\[
\lim_{t\to\infty}\frac{\log(t+1)}{\log(1-F_t(z))}=
\lim_{t\to\infty}\frac{(1-F_t(z))^{1+\alpha}}{(t+1)(1-F_t(z))^\alpha
f \left (F_t(z)\right )}=\frac{\lambda}{2^\alpha a} \,,
\]
and consequently
\[
\lim_{t\to\infty}\frac{d(t)}{\left|1-F_t(z)
\right|^{1+\alpha}\log|1-F_t(z)|}=
\left|\frac{\lambda^{\frac1\alpha}}{2a}\right|\cdot \left|\Im
\left( \frac{\mu}{\alpha\lambda} \right) \right|\,,
\]
which implies assertion (ii).

Let us turn to the case $\beta>\alpha$. As above,
${\beta\in[k\alpha, (k+1)\alpha)}$ for some $k\in\mathbb{N}$. By
Theorem~\ref{th_main_asy}~(iii),
\begin{eqnarray*}
r(z,t)= \frac{\lambda^{\frac{1}{\alpha}}t^{\frac{1}{\alpha}-1}}{2}
\left[\sum_{j=1}^k\left(\begin{array}c 1/{\alpha} \\
j\end{array}\right) \left(h_1(z)\right)^jt^{1-j}  +
\frac{\mu\lambda^{-\frac{\beta}\alpha}t^{1-\frac{\beta}\alpha}}
{\alpha-\beta}
 + r_1(z,t)\right],
\end{eqnarray*}
where $\lim\limits_{t\to\infty}\displaystyle
t^{\frac{\beta}{\alpha}-1} r_1(z,t) =0. $ Substituting this
into~(\ref{dista}) yields
\begin{eqnarray*}
d(t)= \left|\Im\left(\sum_{j=1}^k\left(\begin{array}c 1/{\alpha} \\
j\end{array}\right) \left(h_1(z)\right)^jt^{1-j}  +
\frac{\mu\lambda^{-\frac{\beta}\alpha}t^{1-\frac{\beta}\alpha}}
{\alpha-\beta}  + r_1(z,t) \right)\right| \\
\cdot \frac12\left|\lambda^{\frac1\alpha}\right|
t^{\frac{1}{\alpha}-1}\left|1-F_t(z) \right|^2.
\end{eqnarray*}

For each $z\in\Delta$, there are now two possibilities. One is
that $\Im h_1(z)\not=0$. In this case,
\[
\lim_{t\to\infty}\frac{d(t)}{\left|1-F_t(z) \right|^{1+\alpha}}=
\frac{\left|\lambda^{\frac1\alpha}\right|\left|\Im h_1(z)
\right|}{2\alpha}  \,.
\]
The other possibility is that $h_1(z)$ is real. Then
\begin{eqnarray*}
d(t)= \left|\Im\left( \frac{\mu\lambda^{-\frac{\beta}\alpha}}
{\alpha-\beta}  + t^{\frac{\beta}\alpha-1}r_1(z,t) \right)\right|
\cdot \frac12 \left|\lambda^{\frac1\alpha}\right|
t^{\frac{1-\beta}{\alpha}}\left|1-F_t(z) \right|^2,
\end{eqnarray*}
from which it follows that
\[
\lim_{t\to\infty}\frac{d(t)}{\left|1-F_t(z) \right|^{1+\beta}} =
\frac{\left|\lambda^{\frac\beta\alpha}\right|}{2^\beta(\beta-\alpha)}
\cdot \left|\Im \left(\mu\lambda^{-\frac{\beta}\alpha} \right)
\right|\,.
\]
This implies assertion (iii). \epr

\rema \label{rem1} Theorem~\ref{th_order} shows that the manner of
approaching of different trajectories to their common limit
tangent line essentially depends on the relation between $\alpha$
and $\beta$. For instance, if $\beta<\alpha$ and $\Im
\left(\lambda^{-\frac{\beta}{\alpha}} \mu\right)\not=0$, then by
assertion (i) of Theorem~\ref{th_order}, all the trajectories have
{\it the same} contact order. If $\beta>\alpha$, we see another
phenomenon: by assertion (iii) of Theorem~\ref{th_order}, there
exists a unique trajectory of {\it maximal} contact order. This
has an interesting geometric consequence. The trajectory $\gamma$
which is the pre-image of the real half-axis under $h_1$ has (by
Theorem~\ref{th_order}) contact order at least $\beta$. Suppose
that $\gamma$ is disjoint from the limit tangent line $\ell$.
Since all other trajectories have order $\alpha< \beta$,  each
trajectory starting from a point between $\gamma$ and $\ell$
intersects $\ell$, and for large enough $t$, lies on the opposite
side of $\ell$. \erema

\begin{figure}\centering\label{fig_extremal-line}
    \includegraphics[angle=0,width=5cm,totalheight=4.5cm]{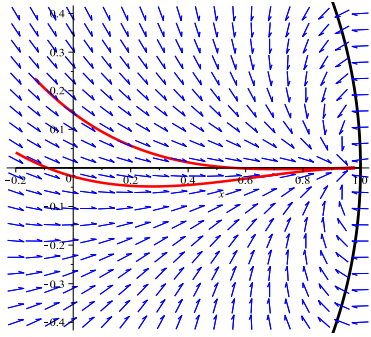}
    \caption{Trajectory of extremal order and trajectory intersecting the tangent line}
\end{figure}

\exa Consider the semigroup $S$ generated by
$f(z)=\displaystyle\frac{(1-z)^2}{4+i(1-z)^2}$\,. Since $a=\frac 1
4$ is real, the limit tangent line $\ell$ coincides with the real
axis. A direct calculation yields $h_1(z)=8+i+iz+\frac{4z}{1-z}$.
Hence, the trajectory defined by
$$\gamma:=\left \{ z \in \Delta: \ \Im
\left(i+iz+\frac{4z}{1-z}\right)=0 \right\}$$ has maximal contact
order, and each trajectory starting from a point between $\gamma$
and the real axis intersects the real axis and eventually lies
below it (see Fig.~\ref{fig_extremal-line}). \eexa

Another implication of Theorem~\ref{th_order} and
Proposition~\ref{cor_asymptotes} relates to the limit curvature of
the semigroup trajectories. Namely, a trajectory
$\gamma_z=\{F_t(z),\ t\ge0 \}\subset\Delta$ has a finite limit
curvature if and only if its image $C\circ\gamma_z$ under the
Cayley transform has an asymptote (cf., \cite{E-S2011}). Moreover,
if that asymptote passes through the point $-1$, then the limit
curvature vanishes. This implies the following.

\begin{corol}\label{cor_asymp_inD}
Let $S=\left\{F_t\right\}_{t\ge0}$ be a semigroup of holomorphic
self-mappings of $\Delta$ generated by $f \in \GD$. Let
$\{\Omega_1, \Omega_2, \Omega_3, \Omega_4, \Omega_5\}$ be the
partition of
 $\Omega=\{(\alpha,\beta): \ 0<\alpha \le 2, \ \beta>0  \}$
as in Proposition~\ref{cor_asymptotes}.

\begin{itemize}
\item[(i)] If  $(\alpha,\beta) \in \Omega_1$, then all the
trajectories of $S$ have null limit curvature.

\item[(ii )] If $(\alpha,\beta) \in \Omega_2$, then each
trajectory has a finite limit curvature (distinct for different
trajectories).

\item[(iii)] If $(\alpha,\beta) \in \Omega_3$, then the trajectory
$\gamma$ defined by the condition $\left. \Im h_1 \right|_\gamma
=0$ is the only trajectory which has a finite limit curvature.
Moreover, this curvature vanishes.

\item[(iv)] If  $(\alpha,\beta) \in \Omega_4$, then all the
trajectories of $\Sigma$ have the same finite limit curvature.
Moreover, this curvature vanishes if and only if $\Im (\mu
\lambda^{-\frac \beta \alpha})=0$.

\item[(v)] If  $(\alpha,\beta) \in \Omega_5$, then there is the
following dichotomy: in the case $\Im (\mu \lambda^{-\frac \beta
\alpha}) =0$, all the trajectories have null limit curvature; in
the case $\Im (\mu \lambda^{-\frac \beta \alpha}) \neq 0$, the
limit curvature is infinite.
\end{itemize}
\end{corol}

If the remainder in the asymptotic expansion of the generator
tends to zero faster than in \eqref{remainders}, the asymptotics
of a semigroup can be estimated more precisely. The next result,
which generalizes \cite[Theorem 1.4(iii)]{E-S2011}, follows from
\eqref{Phi-a}--\eqref{A-B} and transforming
Proposition~\ref{addit_asymp}.

\begin{propo}\label{addit_asymp-d}
Let $S=\left\{F_t\right\}_{t\ge0}$ be a semigroup of holomorphic
self-mappings of $\Delta$ whose generator $f$ has the form
\[
f(z)=a(1-z)^{1+\alpha}+b(1-z)^{1+\alpha+\beta}+R_1(z),
\]
where $R_1\in\Hol(\Delta,\mathbb{C})$ satisfies $\displaystyle
\lim\limits_{z\to1}\frac{R_1(z)}{(1-z)^{2\alpha+1+\varepsilon}}=0$
for some positive $\varepsilon$.
\begin{itemize}
\item[(i)] If $\beta=\alpha$, then there exists a constant $C$
such that for all $z \in \Delta$,
\[
\left(\frac{2}{1-F_t(z)}\right)^{\alpha}= \lambda t + \mu
\log(t+1) + \lambda h(z)+ C + r(z,t),
\]
where $\lim\limits_{t \to 0} r(z,t)=0$.

\item[(ii)] If $\frac\alpha 2<\beta<\alpha$,
 then there exists a
constant $C$ such that for all $z \in \Delta$,
\[
\left( \frac{2}{1-F_t(z)}\right)^{\alpha}=\lambda t+\frac{\mu
\alpha}{\alpha-\beta} (\lambda t)^{1-\frac{\beta}{\alpha}}+\lambda
h(z)+ C + r(z,t),
\]
where $\lim\limits_{t\to\infty}r(z,t)=0$.
\end{itemize}
\end{propo}

\bigskip

\section{Rigidity via order of contact}\label{sect_rigid}
\setcounter{equation}{0}

In this section, we consider two semigroups
$S=\left\{F_t\right\}_{t\ge0}$ and
$S^*=\left\{F^*_t\right\}_{t\ge0}$ acting on $\Delta$. Let $f$ be
the generator of $S$ and $f^*$ the generator of $S^*$. Suppose
that both $f$ and $f^*$ can be represented by~\eqref{circ_4a-e}.
For $z_1,z_2\in\Delta$, let $\ \mathfrak{F}= \left(\{F_t(z_1),\
t\ge0\},\ \{F^*_t(z_2),\ t\ge0\}\right)$ be the pair of semigroups
trajectories. We study the following question: how close can the
trajectories of $S$ and $S^*$ become? Naturally, this question
includes the rigidity problem, i.e., that of determining
conditions which ensure that these semigroups coincide. For this
study we need a modification of Definition~\ref{def_order}.

\begin{defin}\label{def_order_sg}
Let semigroups $S=\left\{F_t\right\}_{t\ge0}$,
$S^*=\left\{F^*_t\right\}_{t\ge0}$ (not necessarily different)
have generators of the form $a(1-z)^{1+\alpha}+ R(z)$ with
$\lim\limits_{t\to \infty}\frac{R(z)}{(1-z)^{1+\alpha}}=0$. Let
$z_1,z_2\in\Delta$. We say that the parameter-related contact
order of $\ \mathfrak{F}= \left(\{F_t(z_1),\ t\ge0\},\
\{F^*_t(z_2),\ t\ge0\}\right)$ is greater than~$\kappa\geq 0$, if
\[
\lim_{t\to\infty} \frac{|F_t(z_1) -F^*_t(z_2)|}
{|1-F_t(z_1)|^{\kappa+1}}=0.
\]

\end{defin}

By Proposition~\ref{th_asy_one_term}, for all $z \in \Delta$, the
limits
\[
\lim_{t\to\infty}
t|1-F_t(z)|^\alpha\quad\mbox{and}\quad\lim_{t\to\infty} t|1- F^*_t
(z)|^\alpha
\]
exist and are finite and nonzero. Therefore, this definition is
symmetric relative to $S$ and $S^*$. Obviously, $|F_t(z_1)
-F^*_t(z_2)|$ is greater than the distance between $F_t(z_1)$ and
the trajectory $\{F^*_s(z_2),s\ge0 \}$. Hence, if the
parameter-related contact order of $\ \mathfrak{F}$ is greater
than $\kappa$, then the contact order between them in a regular
sense (see Definition~\ref{def_order}) is also greater than
$\kappa$.

We consider the two particulary important cases: $S=S^*$ and $z_1
=z_2$. Regarding the case $S=S^*$, it is easy to see from
Corollary~\ref{cor_1} that

if $\beta=\alpha$, then $\displaystyle\lim\limits_{t\to\infty}
\frac{F_t(z)-F_t(0)}{\log(t+1)\left(1-F_t(z)\right)^{1+\alpha}}=0$;

if $\beta<\alpha$, then $\displaystyle\lim\limits_{t\to\infty}
\frac{F_t(z)-F_t(0)}{\left(1-F_t(z)\right)^{1+\beta}}=0$;

if $\beta>\alpha$, then $\lim\limits_{t\to\infty}
t^{\frac{\beta}{\alpha}-1}\left(\displaystyle\frac{F_t(z)-F_t(0)}{\left(1-F_t(z)\right)^{1+\alpha}}
- \displaystyle\frac{\lambda^\frac{1}{\alpha}
\sigma(w)}{\alpha}\right)=0$.

By Definition~\ref{def_order_sg}, this implies the following fact.

\begin{propo}\label{similar_inD}
Let $S=\{F_{t}\}_{t\geq 0}\subset \Hol(\Delta)$ be a semigroup
generated by a mapping $f \in \GD$, and let $z_1,z_2\in\Delta$.
\begin{itemize}
\item[(i)] If $\beta<\alpha$, then the parameter-related contact
order of $\ \mathfrak{F}$ is greater than $\beta$.

\item[(ii)] If $\beta\ge\alpha$, then for any $\varepsilon>0$, the
parameter-related contact order of $\ \mathfrak{F}$ is greater
than $\alpha-\varepsilon$.
\end{itemize}
\end{propo}

\rema \label{rem2} Just from the triangle inequality, it follows
that the parameter-related contact order of each pair of
trajectories of the same semigroup cannot be less than the contact
order of any one of them with the limit tangent line. Comparing
Proposition~\ref{similar_inD} with Theorem~\ref{th_order}, we see
that if $0<\beta<\alpha$ and $\Im
\left(\lambda^{-\frac{\beta}{\alpha}}\mu\right)\not=0$, then the
parameter-related contact order of two trajectories is actually
greater than the contact order of any of them with the limit
tangent line. Roughly speaking, this means that each trajectory is
closer to all other trajectories than to the tangent. This can be
the case only when all trajectories approach their common limit
tangent line from the same side. \erema

\begin{figure}\centering\label{one-side}

    \includegraphics[angle=0,width=5cm,totalheight=4.5cm]{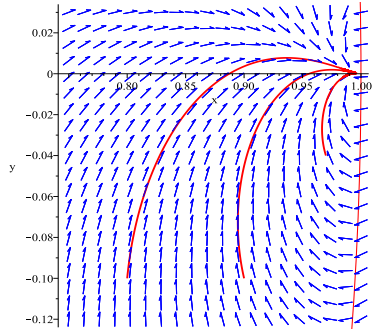}
    \caption{All trajectories approach the tangent line from the same side}
\end{figure}

\exa Consider $f\in\GD$ defined by $f(z):=(1-z)^2 + i(1-z)^{2.5}$,
i.e., $\alpha=1,\ \beta=0.5<\alpha,\ a=1$ and $b=i$. Since $\arg
a=0$, the limit tangent line coincides with the real axis. In
addition, $\lambda=2$ and $\mu=i\sqrt{2}$, so that $\Im
\left(\lambda^{-\frac{\beta}{\alpha}}\mu\right)=1\not=0$.
Fig.~\ref{one-side} shows the direction field in the part of
$\Delta$ bounded by $0.75<\Re z<1$ and $-0.12<\Im z<0.03$. All
trajectories approach the real axis from the upper half-plane.
\eexa

We now turn to the case $z_1=z_2$. We are interested in applying
contact order to the rigidity problem.

\begin{theorem}\label{th_weak_rig}
Let $S$ be a semigroup generated by a mapping $f\in\GD$ with
$\beta\le\alpha$, and let $S^*$ be a semigroup generated by
$f^*(z)=f(z)+c(1-z)^{1+\alpha+\beta}$. If for some $z \in \Delta$,
the parameter-related contact order of $\ \mathfrak{F}=
\left(\{F_t(z),\ t\ge0\},\ \{F^*_t(z),\ t\ge0\}\right)$ is greater
than $\beta$, then $c=0$; so the semigroups coincide.
\end{theorem}

\pr It follows from Proposition~\ref{th_asy_one_term} that
\begin{equation}\label{auxil1}
\lim_{t\to\infty} t\left|1-F_t(z)\right|^\alpha =
\lim_{t\to\infty} t \left|1- F^*_t(z)\right|^\alpha
=\frac{2^\alpha}\lambda\,.
\end{equation}
Let $\beta<\alpha$. Consider the quotient
\begin{eqnarray*}
&&\frac{F_t(z)- F^*_t(z)}{(1-F_t(z))^{1+\beta}} = \\
=&& t^{\frac{\beta-1}\alpha}\left(\frac1{1-F_t(z)}- \frac1{1-
F^*_t(z)}\right)\cdot \frac{1- F^*_t(z)}{1-F_t(z)}\cdot \frac1
{(1-F_t(z))^{\beta-1}t^{\frac{\beta-1}\alpha}} \,.
\end{eqnarray*}
Formula~\eqref{auxil1} implies that the last two factors have
finite nonzero limits. In addition, by Theorem~\ref{th_main_asy}
(i),
\[
\frac1{1-F_t(z)}=\frac{(\lambda t)^{\frac{1}{\alpha}}}{2}
\left(1+\frac{\mu} {\alpha-\beta}(\lambda
t)^{-\frac{\beta}{\alpha}} + r_1(z,t)\right)
\]
and
\[
\frac1{1- F^*_t(z)}=\frac{(\lambda t)^{\frac{1}{\alpha}}} {2}
\left(1+\frac{\mu^*} {\alpha-\beta}(\lambda
t)^{-\frac{\beta}{\alpha}} + r^*_1(z,t)\right),
\]
where $\mu^* =\frac{2^\beta(b+c)}a$ and $\lim\limits_ {t\to\infty}
\displaystyle t^{\frac{\beta}{\alpha}} r_1(z,t) =
\lim\limits_{t\to\infty}\displaystyle t^{\frac{\beta}{\alpha}}
r^*_1(z,t) =0.$ Therefore,
\begin{eqnarray}\label{aux2}
&&\frac1{1-F_t(z)}-  \frac1{1- F^*_t(z)} \nonumber \\
=&& \frac{(\lambda t)^{\frac{1-\beta}{\alpha}}}{2} \left(\frac{\mu
- \mu^*} {\alpha-\beta} + (\lambda t)^{\frac{\beta}{\alpha}}
\left(r_1(z,t) - r^*_1(z,t)\right)\right).
\end{eqnarray}
Thus, if the parameter-related contact order of the pair
$\mathfrak{F}$ is greater than $\beta$, then $\mu-\mu^*=0$, and
the assertion follows.

The case $\beta=\alpha$ can be treated analogously using
assertion~(ii) of Theorem~\ref{th_main_asy}. \epr

Our next result concerns the rigidity problem  in the case
$\beta>\frac {\alpha}{2}$. We prove the coincidence of semigroups
under an essentially weaker local condition.

\begin{theorem}\label{th_strong_rig}
Let semigroups $S$ and $S^*$ be generated, respectively, by
\[
f(z)=a(1-z)^{1+\alpha}+b(1-z)^{1+\alpha+\beta}+R_1(z)
\]
and
\[
f^*(z)=a(1-z)^{1+\alpha}+ b^*(1-z)^{1+\alpha+\beta}+ R^*_1(z)
\]
with $a\not=0$ and $R_1, R^*_1\in\Hol(\Delta,\mathbb{C})$. Suppose
that either

(i) $\frac\alpha2<\beta\le\alpha$ and $\displaystyle
\lim\limits_{z\to1}\frac{R_1(z)}{(1-z)^{1+2\alpha+\varepsilon}}=
\lim\limits_{z\to1}
\frac{R^*_1(z)}{(1-z)^{1+2\alpha+\varepsilon}}=0$ for some
positive $\varepsilon$, or

(ii) $\beta>\alpha$ and $\displaystyle
\lim\limits_{z\to1}\frac{R_1(z)}{(1-z)^{1+\alpha+\beta}}=
\lim\limits_{z\to1}
\frac{R^*_1(z)}{(1-z)^{1+\alpha+\beta}}=0$.\\
If there exist $\theta\in[0,2\pi]$ and an open set $\mathcal{U}
\subset \Delta$ such that
\begin{equation}\label{str_rig1}
\lim_{t\to\infty} t^{1+\frac1\alpha}\Re e^{i\theta} \left(F_t(z) -
F^*_t(z) \right)=0
\end{equation}for all $z \in \mathcal{U}$,
then the two semigroups coincide.
\end{theorem}

\pr  Consider the case $\frac\alpha2<\beta<\alpha$. By our
assumptions, both semigroups satisfy the asymptotic
expansion~\eqref{as_less}. Since
\[
\frac1{1-F_t(z)} - \frac1{1- F^*_t(z)} =\frac{F_t(z)- F^*_t(z)}
{(1-F_t(z))(1-F^*_t(z))},
\]
the parameter-related contact order of $\mathfrak{F}=
\left(\{F_t(z),\ t\ge0\},\ \{F^*_t(z),\ t\ge0\}\right)$ is
positive by formula \eqref{aux2}. Furthermore,
\begin{eqnarray*}
&&\left( \frac{2}{1-F_t(z)}\right)^{\alpha} - \left(
\frac{2}{1- F^*_t(z)}\right)^{\alpha} \\
&=& \left( \frac{2}{1-F_t(z)}\right)^{\alpha} \cdot \frac{F^*_t(z)
- F_t(z)} {1- F^*_t(z)} \cdot \frac{1- \left( 1+ \frac{F^*_t(z) -
F_t(z)} {1-F^*_t(z)} \right)^\alpha } { \frac{F^*_t(z) - F_t(z)}
{1- F^*_t(z)}},
\end{eqnarray*}
where the last factor tends to $-\alpha$ as $t\to\infty$. On the
other hand, by Proposition~\ref{addit_asymp-d},
\begin{eqnarray*}
&&\left( \frac{2}{1-F_t(z)}\right)^{\alpha} - \left(
\frac{2}{1- F^*_t(z)}\right)^{\alpha} \\
&& = \frac{(\mu- \mu^*) \alpha}{\alpha-\beta} (\lambda
t)^{1-\frac{\beta}{\alpha}}+\lambda \left(h(z)- h^*(z)\right) + C
- C^* + r(z,t) - r^*(z,t),
\end{eqnarray*}
where $C, C^*$ are constants and $h, h^*$ are the K{\oe}nigs
functions for $S$ and $S^*$, respectively.

Assume condition~\eqref{str_rig1}. Combining the last two
displayed formulas with \eqref{repr_sh}, we conclude that
\[
\lim_{t\to\infty} \Re \left(
\frac{e^{i\theta}}{\lambda^{1+\frac1\alpha}} \left(\frac{(\mu-
\mu^*) \alpha}{\alpha-\beta} (\lambda
t)^{1-\frac{\beta}{\alpha}}+\lambda \left(h(z)- h^*(z)\right) + C
- C^* \right)\right)=0.
\]
This is possible only if the coefficient of
$t^{1-\frac{\beta}{\alpha}}$ vanishes, in which case
\[
\Re \left( \frac{e^{i\theta}}{\lambda^{1+\frac1\alpha}}
\left(\lambda \left(h(z)- h^*(z)\right) + C - C^*
\right)\right)=0.
\]
Therefore, the function ${h(z)- h^*(z)}$ is constant. Since
$h(0)=h^*(0)=0,$ we get $h(z)= h^*(z)$. Now by (\ref{h}), we
conclude that $f\equiv f^*$.

The case $\beta>\alpha$ can be treated similarly. \epr

Suppose now that \eqref{str_rig1} holds for all
$\theta\in[0,2\pi]$. By Proposition~\ref{th_asy_one_term},
$\displaystyle t^{1+\frac1\alpha}\thicksim
\frac1{(1-F_t(z))^{\alpha+1}}$. Thus, we get the following
consequence.

\begin{corol}
If under conditions of Theorem~\ref{th_strong_rig}, the
parameter-related contact order of $\ \mathfrak{F}$ is greater
than $\min(\alpha,\beta)$, then the semigroups $S$ and $S^*$
coincide.
\end{corol}

Let return to the formulation of Theorem~\ref{th_strong_rig}. It
seems that the requirement on the remainders in assertion~(i) is
too strong and should be replaced by $\displaystyle
\lim\limits_{z\to1}\frac{R_1(z)}{(1-z)^{1+\alpha+\beta}}=
\lim\limits_{z\to1} \frac{R^*_1(z)}{(1-z)^{1+\alpha+\beta}}=0$, as
in assertion (ii). Moreover, the rigidity condition
\eqref{str_rig1} does not include $\beta$ at all. These
considerations lead to the following natural conjecture.

\begin{conj}  Let $S$ and $S^*$ be semigroups generated by functions
of the form $a(1-z)^{1+\alpha}+ R(z)$. If for some $\varepsilon>0$
the remainders are $O((1-z)^{1+\alpha+\varepsilon})$, then
condition \eqref{str_rig1} implies the coincidence of the
semigroups.
\end{conj}

\bigskip

\section{Appendix}
\setcounter{equation}{0}

We complete our analysis with assertions which give more
information about the asymptotic behavior of semigroups but are
different in nature.

\begin{propo}\label{propo_add}
Let $\Sigma=\{\Phi_{t}\}_{t\geq 0}\subset\Hol(\Pi)$ be a semigroup
generated by $\phi \in \GP$. Then
\[
\lim_{t\to\infty} (t+1)^{\frac{\beta}{\alpha}}
\left((\Phi_t(w)+1)^{\alpha}-(\Phi_t(1)+1)^{\alpha} -
\lambda\sigma(w)\right)
=\mu\lambda^{1-\frac{\beta}{\alpha}}\sigma(w).
\]
\end{propo}

\pr We just calculate the limit:
\begin{eqnarray*}
&&\lim_{t\to\infty}
(t+1)^{\frac{\beta}{\alpha}}\left((\Phi_t+1)^{\alpha}(w)-(\Phi_t+1)^{\alpha}(1)
-\lambda\sigma(w)\right) \\
&= &\lim_{t\to\infty} (t+1)^{\frac{\beta}{\alpha}}\cdot \int_1^w
\left((\Phi_t(z)+1)^{\alpha}-\lambda\sigma(z)\right)'dz \\
&=& \lim_{t\to\infty} (t+1)^{\frac{\beta}{\alpha}}\cdot \int_1^w
\frac{\alpha(\Phi_t(z)+1)^{\alpha-1}\phi\left(\Phi_t(z)\right)-\lambda}{\phi(z)}
dz.
\end{eqnarray*}

Since
\begin{eqnarray*}
&&\lim_{t\to\infty}(t+1)^{\frac{\beta}{\alpha}}
\left(\alpha(\Phi_t(z)+1)^{\alpha-1}\phi\left(\Phi_t(z)\right)-\lambda\right)
\\
&=&\lim_{t\to\infty}\left(\frac{(t+1)^{\frac{1}{\alpha}}}{(\Phi_t(z)+1)}\right)^{\beta}
\left(\alpha
B+\frac{\rho_{1}((\Phi_t(z))}{(\Phi_t(z)+1)^{1-\alpha-\beta}}\right)
=\mu\lambda^{1-\frac{\beta}{\alpha}},
\end{eqnarray*}
we conclude that
\[
\lim_{t\to\infty}
(t+1)^{\frac{\beta}{\alpha}}\left((\Phi_t(w)+1)^{\alpha}-(\Phi_t(1)+1)^{\alpha}
-
\lambda\sigma(w)\right)=\mu\lambda^{1-\frac{\beta}{\alpha}}\sigma(w),
\]
which completes the proof. \epr

The particular case $\alpha=\beta=1$ is contained in \cite[Theorem
4.1(ii)]{E-S2011}. Transferring, as above,
Proposition~\ref{propo_add} to semigroups acting in $\Delta$
yields the following result.

\begin{corol}
Let $S=\left\{F_t\right\}_{t\ge0}$ be a semigroup of holomorphic
self-mappings of $\Delta$ generated by $f \in \GD$. Then
\[
\lim_{t\to\infty} (t+1)^{\frac{\beta}{\alpha}} \left( \frac 1
{(1-F_t(z))^\alpha } - \frac 1 {(1-F_t(0))^\alpha } -
\frac{\lambda h(z)}{2^\alpha}\right)
=\frac{\mu\lambda^{1-\frac{\beta}{\alpha}} h(z)}{2^\alpha}\,.
\]
\end{corol}

\noindent{\bf Acknowledgment.} The study of semigroups whose
generators have a non-integer power asymptotic expansion was first
proposed by David Shoikhet in joint papers \cite{E-S-Y} and
\cite{E-K-R-S}. The authors are grateful to him for fruitful
collaboration.



\end{document}